\documentclass[%
 aip,
 amsmath,amssymb,
 groupedaddress,
 reprint,%
onecolumn
]{revtex4-1}

\usepackage{graphicx}
\usepackage{dcolumn}
\usepackage{bm}

\usepackage[utf8]{inputenc}
\usepackage[T1]{fontenc}
\usepackage{mathptmx}
\usepackage{etoolbox}

\usepackage[english]{babel}
\usepackage[utf8]{inputenc}
\usepackage[T1]{fontenc}
\usepackage{amsmath,amssymb}
\usepackage{graphicx}
\usepackage{float}

\usepackage{transparent} 
\usepackage{eso-pic} 
\usepackage{tikz} 
\usetikzlibrary{}
\usepackage{caption} 
\usepackage[labelformat=simple]{subcaption}

\usepackage{ragged2e} 
\DeclareCaptionJustification{justified}{\justifying}

\usepackage{xcolor} 

\usepackage[export]{adjustbox}  
\usepackage{fancyhdr}  
\usepackage{hyperref}  
\usepackage{datetime}  
\usepackage{array}  

\usepackage[symbol]{footmisc}

\makeatletter
\def\@email#1#2{%
 \endgroup
 \patchcmd{\titleblock@produce}
  {\frontmatter@RRAPformat}
  {\frontmatter@RRAPformat{\produce@RRAP{*#1\href{mailto:#2}{#2}}}\frontmatter@RRAPformat}
  {}{}
}%
\makeatother

\graphicspath{{./FIGURES/}} 

\begin{document}

\preprint{AIP/123-QED}

\title{Data-Driven Nonlinear Model Reduction to Spectral Submanifolds \\via Oblique Projection}


\author{Leonardo Bettini, Bálint Kaszás, Bernhard Zybach, Jürg Dual}
\author{George Haller}
 \altaffiliation[]{Author to whom correspondance should be addressed: georgehaller@ethz.ch}
\affiliation{Institute for Mechanical Systems, ETH Z\"urich, Leonhardstrasse 21, 8092 Z\"urich, Switzerland}

%
%

\begin{abstract}
\begin{center}
\Large{(To appear in \textbf{Chaos})}\\
\end{center}

The dynamics in a primary Spectral Submanifold (SSM) constructed over the slowest modes of a dynamical system provide an ideal reduced-order model for nearby trajectories. Modeling the dynamics of trajectories further away from the primary SSM, however, is difficult if the linear part of the system exhibits strong non-normal behavior. Such non-normality implies that simply projecting trajectories onto SSMs along directions normal to the slow linear modes will not pair those trajectories correctly with their reduced counterparts on the SSMs. In principle, a well-defined nonlinear projection along a stable invariant foliation exists and would exactly match the full dynamics to the SSM-reduced dynamics. This foliation, however, cannot realistically be constructed from practically feasible amounts and distributions of experimental data. Here we develop an oblique projection technique that is able to approximate this foliation efficiently, even from a single experimental trajectory of a significantly non-normal and nonlinear beam.

\end{abstract}

\maketitle

\begin{quotation}

Data-driven reduced-order models seeking to describe, predict and control complex physical systems attract growing interest in various areas of applied science and engineering. Reduction to spectral submanifolds (or SSMs) provides a mathematically justified model reduction approach that captures the nonlinear dynamics on low-dimensional, attracting invariant manifolds in the phase space of the system. Nearby trajectories approach SSMs and quickly synchronize with their internal dynamics, which then serve as accurate reduced models for the longer-term dynamics of the full system. Therefore, establishing a proper correspondence between a trajectory outside an SSM and its target trajectory within the SSM is essential. This task becomes challenging when the linear part of the system is strongly non-normal, necessitating an oblique projection of trajectories onto their targets within the SSM. Here, we develop a data-driven identification of this oblique projection and allows accurate SSM reduction from a few experimentally observed trajectories of highly non-normal and nonlinear systems.

\end{quotation}

\section{Introduction}
Reduced-order models have become crucially important in studying high-dimensional nonlinear dynamical systems across various fields of applied science and engineering. They allow us to reduce computational cost while retaining physical interpretability and enabling fast parametric analysis, especially when only experimental data are available. Numerous model reduction methods have been proposed using vastly different approaches (see \citet{benner2015,rowley2017,taira2017,brunton2020,touze2021,ghadami2022} for recent reviews). 

Proper orthogonal decomposition (POD), followed by Galerkin projection is a popular and simple approach to reduce the dimension of nonlinear dynamical systems. It employs a normal projection onto the most energetic modes of the system. However, the main limitation of projecting onto linear subspaces is that those subspaces are not invariant under the nonlinear dynamics (\citet{Ohlberger_2016}). This motivates different techniques that project the dynamics onto general invariant manifolds, whose internal dynamics serve as a reduced order model for the system. We employ the theory of spectral submanifolds (SSMs), which allows us to perform accurate model reduction for essentially nonlinear phenomena. Such SSMs rigorously extend the concept of nonlinear normal modes (NNMs) originally introduced by \citet{shaw_pierre_93,shaw_pierre_94} and \citet{shaw_pierre_99}. They have been mathematically formalized by \citet{cabre2005} and \citet{haller2016} for systems with asymptotically stable stationary states. In particular, the primary SSM is the smoothest continuation of a spectral subspace of the linearized system at the stationary state under the inclusion of nonlinearities. Its existence, uniqueness and persistence in both autonomous and non-autonomous systems depend on certain non-resonance conditions among eigenvalues of the spectrum of the linearization (see \citet{haller2016}). 

The internal dynamics within the primary SSM tangent to the spectral subspace spanned by the slowest modes is an ideal nonlinear reduced model for the system. Indeed, such an SSM attracts all nearby trajectories, each of which then synchronizes exponentially fast with a particular trajectory within the SSM. Therefore, properly identifying such a target trajectory on the SSM for each off-SSM initial condition is crucial for effective reduced-order modelling.

Classic linearization results near hyperbolic fixed points guarantee that off-SSM initial conditions converging the fastest to a given on-SSM initial condition form an $(n-d)$-dimensional surface, where $n$ is the dimension of the full system and $d$ is the dimension of the SSM (see \citet{haller_2024}). This surface is usually called the stable fiber associated with the on-SSM initial condition, which is referred to as the base point of the stable fiber. Stable fibers are known to be as smooth as the dynamical system and provide a smooth foliation of a neighborhood of the SSM. This stable foliation is invariant as fibers are mapped into fibers by the full flow map (see \citet{szalai_2023} for more discussion).

If the initial condition lies close to the SSM, the SSM is nearly flat (this is often the case in a delay-embedding setting, see \citet{cenedese2022}) and the linear part is close to normal, then projecting such initial condition orthogonally onto the SSM yields an accurate reduced initial condition whose trajectory within the SSM serves as a target for the trajectory of the full system. Such a normal projection approximates the stable fibers of the SSM as planes aligned with the direction of the fast subspace. In our work, \textit{orthogonality} or \textit{normality} refers to the mutual orientation of eigenspaces in the phase space. The above SSM-reduction combined with orthogonal projection has proven successful in both equation- and data-driven problems, representing essentially nonlinear features. In the data-driven case, unforced, decaying data is employed to fit the parametrization of the primary SSM and its reduced dynamics. External forcing can then be taken into account to predict the nonlinear response of the system (see \citet{haller2016,cenedese2022,cenedese2022_mechanical,joar_2023,mattia_toappear,haller_2024}), even when non-smoothness is involved (\citet{bettini_2024}). 

Recent work by \citet{haller2023_fractional} has revealed the existence of an additional, infinite family of fractional (or secondary) SSMs in $C^\infty$ dynamical systems. These are tangent to the same spectral subspace as the primary SSM, but they fill an entire open set of the domain of attraction of the fixed point. Fractional SSMs are of lower smoothness class than the primary SSM, which is indicated by the appearance of non-integer powers in their parametrization. This new class of SSMs arises in modeling transitions between isolated states (\citet{balint_2024}) or in experimental problems, as general initial conditions lie on a fractional SSM with probability one.  \citet{haller2023_fractional} also extend the theory of SSMs to spectral subspaces with both stable and unstable directions (mixed-mode SSMs).

Among other techniques exploiting manifold reduction through normal projection, \citet{lee_2019,Champion_2019,Fresca_2021,Conti_2023} and \citet{romor_2023} employ an autoencoder. Specifically, the reduced dynamics and manifold parametrization are learned in the process of encoding the data into a lower-dimensional space and subsequently decoded in the original space (\citet{Goodfellow_2016}). As an alternative, \citet{geelen_2023} and \citet{benner_2024} orthogonally project onto a manifold expressed as a graph over a slow subspace computed from proper orthogonal decomposition. 

However, when the non-normality of the system becomes more significant and the slow manifold is far from being flat, a more refined projection of initial conditions to the base point of the stable fiber containing them becomes essential. Here, we call a system \textit{non-normal} if the slow subspace, containing the modes of the linearized dynamics retained in the reduced-order model, is not orthogonal to the fast subspace, which collects the remaining faster modes. Such behavior is common in fluid dynamics, as described by \citet{Trefethen_1993,Schmid_2002,Trefethen_2005}. In structural vibrations, experimental data from resonance-based gravity measurements (\citet{brack2022}) and hydrogel vibrations (\citet{Yerrapragada_2024}), also show clear signs of non-normal behavior.

In these problems, a more refined projection onto the SSM is required along the stable fibers emanating from the SSM. Fiber reconstruction from data relies heavily on the exponentially decaying transients, which make stable fiber reconstruction challenging from experiments. Indeed, as the examples of \citet{szalai_2020,szalai_2023} show, an accurate reconstruction of the stable foliation near an SSM requires data amounts and densities that are unrealistic to obtain in a real experimental setting. For instance, \citet{szalai_2023} approximates the local stable foliation of a 1D SSM in a 2D system using 500 trajectories whose initial conditions are uniformly distributed in the phase space. Neither the number nor the placement of these initial conditions is feasible in an experimental setting. Also, the numerical burden of constructing stable fibers from data increases significantly with the dimension of the problem. Indeed, while SSM-reduction deals with a single $d$-dimensional manifold, stable fiber construction targets a $d$-dimensional family of $(n-d)$-dimensional manifolds. In practical problems, one actually deals with continuum, for which we have $n = \infty$.

In the context of input-output systems, the stable fibers can be approximated through oblique projections directly from the equations. Given a target subspace spanned by the reduced coordinates of the envisioned reduced-order model, oblique projections allow one to connect the portion of the phase space outside that subspace to base points contained in the subspace. Balanced truncation (\citet{Mullis_1976,Moore_1981}) and balanced POD (\citet{ROWLEY_2005,Bagheri_2009,BARBAGALLO_2009}) provide a linear approximation, while \citet{Scherpen_1993,Scherpen_1994} extended the balancing idea between observable and controllable states to nonlinear projections. However, in this work we focus on data-driven reduced-order modeling: for a more extensive review of equation-driven balancing techniques, refer to \citet{Gugercin_2004,Benner_breiten_2017}. 

Recent work by \citet{otto_2023} constructs nonlinear oblique projections directly from data via an autoencoder neural network. The architecture of the encoder is constrained to represent a projection, while the decoder parametrizes the manifold over which the projection takes place. The direction of the projection is then learned by minimizing a loss function that measures the difference between the time derivative of the projected trajectories and the reduced dynamics on the manifold, evaluated at the projected points. However, this method also requires large amounts of data, such as $1000$ training trajectories in an example given by \citet{otto_2023}. In some cases, the direction of stable fibers is primarily influenced by non-normality of the slow and fast eigenspaces, rather than by the presence of nonlinearities. In such situation, a linear oblique projection is sufficient for achieving effective reduced-order models (see \citet{BARBAGALLO_2009,AHUJA_2010,Illingworth_2011,Benner_2018,Otto_2022,Otto_2023_covariance}).

Another form of simplification arises when the eigenvalues of the slow and fast subspaces are known. \citet{axas_2023} show this to be the case for delay embeddings of systems with a known linear spectrum. Optimizing the number and magnitude of delay used, one can then achieve near-normality in the delay-embedding space, and use normal projection to the slow subspace with high accuracy.

In the present work, we employ linear, oblique projections onto primary SSMs to construct SSM-reduced models for nonlinear and non-normal dynamical systems. The method requires a low number of trajectories and hence applies to both numerical and experimental data, as we will show. The idea is to identify the overall impact of the fast dynamics in the data, rather than identify the fast subspace itself. Specifically, the influence of the fast dynamics on a decaying trajectory manifests itself in oscillations of the backbone curve, i.e. the curve of instantaneous amplitudes as functions of instantaneous frequencies in an oscillatory decaying signal.
Based on this observation, we construct an oblique projection along stable fibers by finding the linear mapping from the full phase space to the slow spectral subspace that minimizes the oscillations in the backbone curve of a decaying trajectory observed under that projection. This procedure allows us to effectively approximate the overall direction of the fibers even from a \textit{single} decaying trajectory, as opposed to the vast amounts and densities required for the data by other methods we have surveyed above.

The structure of this paper is as follows. In section \ref{motivation} we first provide an explanation for oscillatory backbone curves, which are present already in non-normal linear systems. We identify the two key contributors to this phenomenon: fractional SSMs and non-normal eigenspaces. Section \ref{method} outlines our model reduction procedure based on obliquely projected reduced coordinates over which the primary SSM is constructed. Finally, in section \ref{results}, we apply our oblique SSM reduction method to both numerical and experimental data sets.

\section{Motivation}\label{motivation}
Let us consider the nonlinear dynamical system
\begin{equation}
    \label{eq:basic_eq}
    \dot{\mathbf{x}} = \mathbf{A}\mathbf{x} + \mathbf{f}_\mathrm{nl}(\mathbf{x}), \qquad \mathbf{x}\in \mathbb{R}^n,
\end{equation}
with a fixed point at $\mathbf{x}=0$, linear part $\mathbf{A}\in \mathbb{R}^{n\times n}$ and nonlinearities $\mathbf{f}_\mathrm{nl} =  \mathcal{O}(|\mathbf{x}|^2)$. We assume that the fixed point is linearly asymptotically stable, i.e., $\text{Re } \lambda_i < 0 $ for the eigenvalues $\lambda_i$ of $\mathbf{A}$ for $i=1,...,n$. Moreover, we assume that $\mathbf{A}$ is semisimple and that its eigenvalues satisfy the nonresonance condition 
\begin{equation}
\lambda_j \neq \sum_{k = 1}^n m_k \lambda_k,\quad m_k \in \mathbb{N}, \quad \sum_{k=1}^n m_k \geq 2, \quad j = 1, \dots, n.
\end{equation}
SSM-based model reduction constructs a slow attracting SSM $\mathcal{W}\left( E\right)$ of the fixed point and restricts the full dynamics to this low-dimensional manifold. The reduced coordinates on the SSM are defined as a projection of the coordinates of the phase space onto the slow spectral subspace $E$ of $\mathbf{A}$ that is tangent to $\mathcal{W}\left( E\right)$ at the fixed point. The internal dynamics of the slow SSM are then represented as the dynamics of the reduced coordinates. Normal projection of coordinates onto the slowest spectral subspace is justified when the spectral subspace $F$ spanned by the remaining faster-decaying eigenmodes of $\mathbf{A}$ is normal to the slow spectral subspace tangent to the SSM. When this is not the case, an oblique projection along $F$ is needed to account for the effects due to non-normality.

As an illustration, consider the linear system
\begin{equation}\label{2d_example}
	\mathbf{\dot{x}} = \mathbf{A}\mathbf{x},\quad \mathbf{x} = (x_1, x_2,)^\mathrm{T},\quad \mathbf{A} = 
	\begin{pmatrix}
		-\alpha & 0 \\
		0 & -\beta\\ 
	\end{pmatrix}, \quad\beta > \alpha > 0, 
\end{equation}
where the slow and fast subspaces $E$ and $F$ are 1D and are spanned by the eigenvectors corresponding to the eigenvalues $\lambda_1 = -\alpha$ and  $\lambda_2 = -\beta$:
\begin{equation}
\mathbf{v}^{(1)} = 
\begin{pmatrix}
	1 \\ 0
\end{pmatrix}, \quad
\mathbf{v}^{(2)} = 
\begin{pmatrix}
	0 \\ 1
\end{pmatrix}, 
\quad E = \text{span}\left\{ \mathbf{v}^{(1)}\right\}, \quad F = \text{span}\left\{ \mathbf{v}^{(2)}\right\}.
\end{equation}\label{2d_eigenvectors}
Note that the linear system is normal, i.e. $\mathbf{v}^{(1),\mathrm{T}}\cdot \mathbf{v}^{(2)} =  0$. Moreover, the reduced dynamics of fractional SSMs filling the phase space are identical to the reduced dynamics of the primary SSM, $\mathcal{W}\left( E\right) \equiv E$. Consequently, even though a general initial condition lies on a fractional SSM, the normal projection of the trajectory that it generates is accurately modeled by the internal dynamics of $\mathcal{W}\left( E\right)$
\begin{equation}\label{internal_SSM_2d}
	x_1(t) = x_{10}e^{-\alpha t}.
\end{equation}

\begin{figure*}
    \centering
    \subfloat[\label{1d_ex_orth}]{
        \includegraphics[scale=0.35]{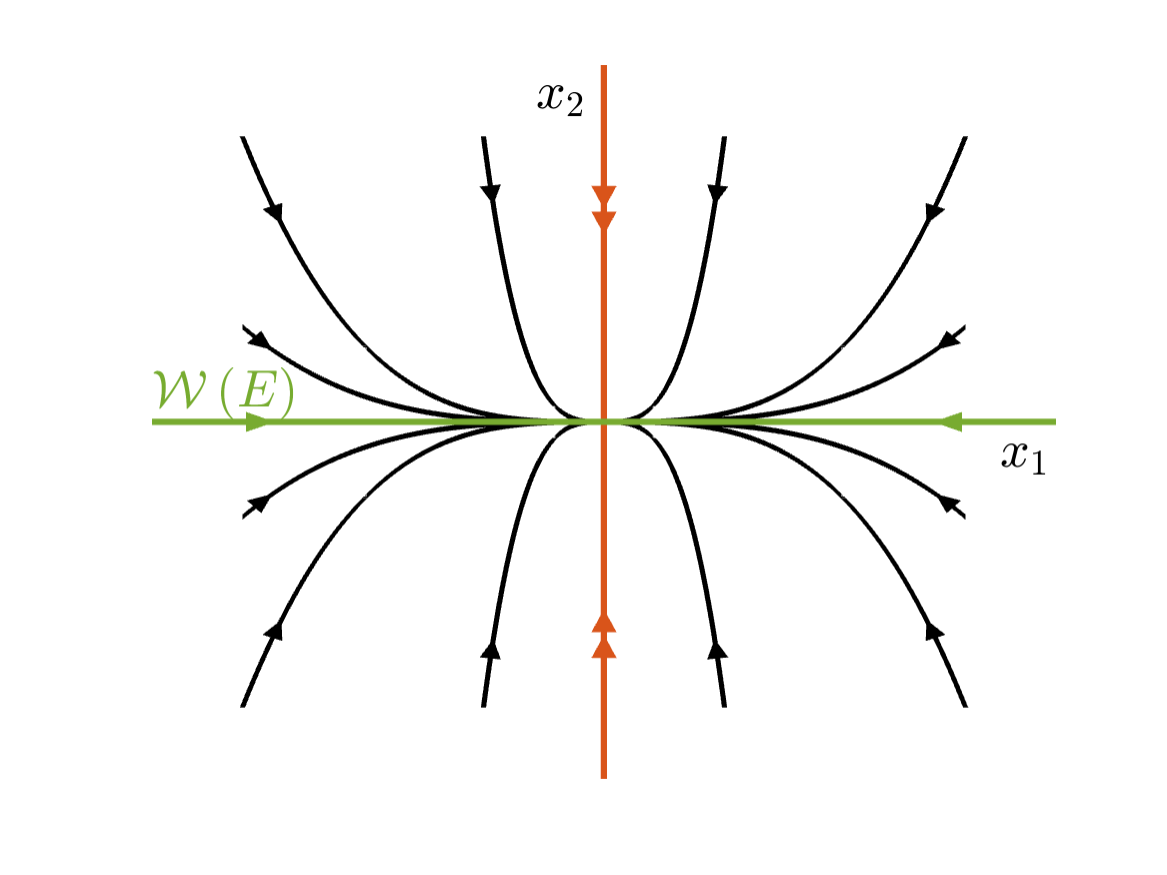}
    }
    \quad
    \subfloat[\label{1d_ex_obl}]{
        \includegraphics[scale=0.35]{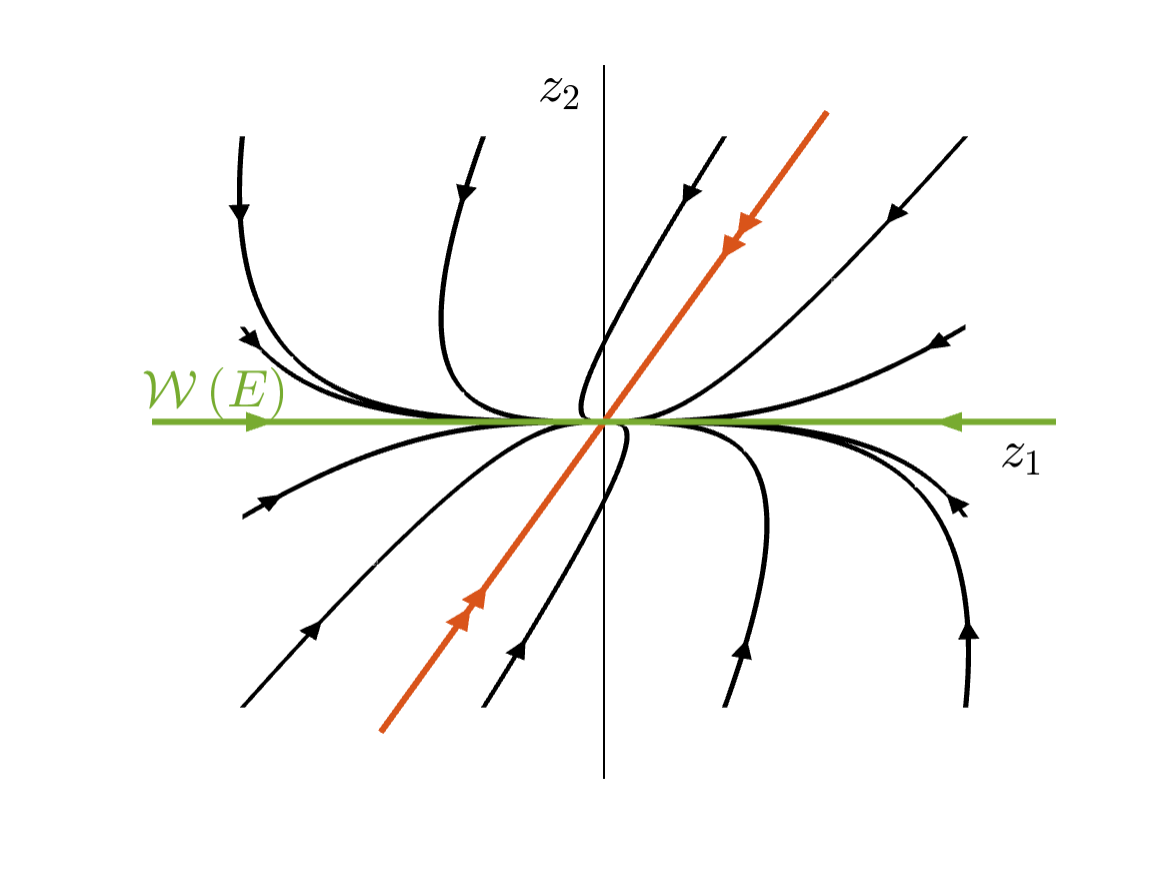}
    }
    \quad
    \subfloat[\label{2d_ex_backones}]{
        \includegraphics[scale=0.35]{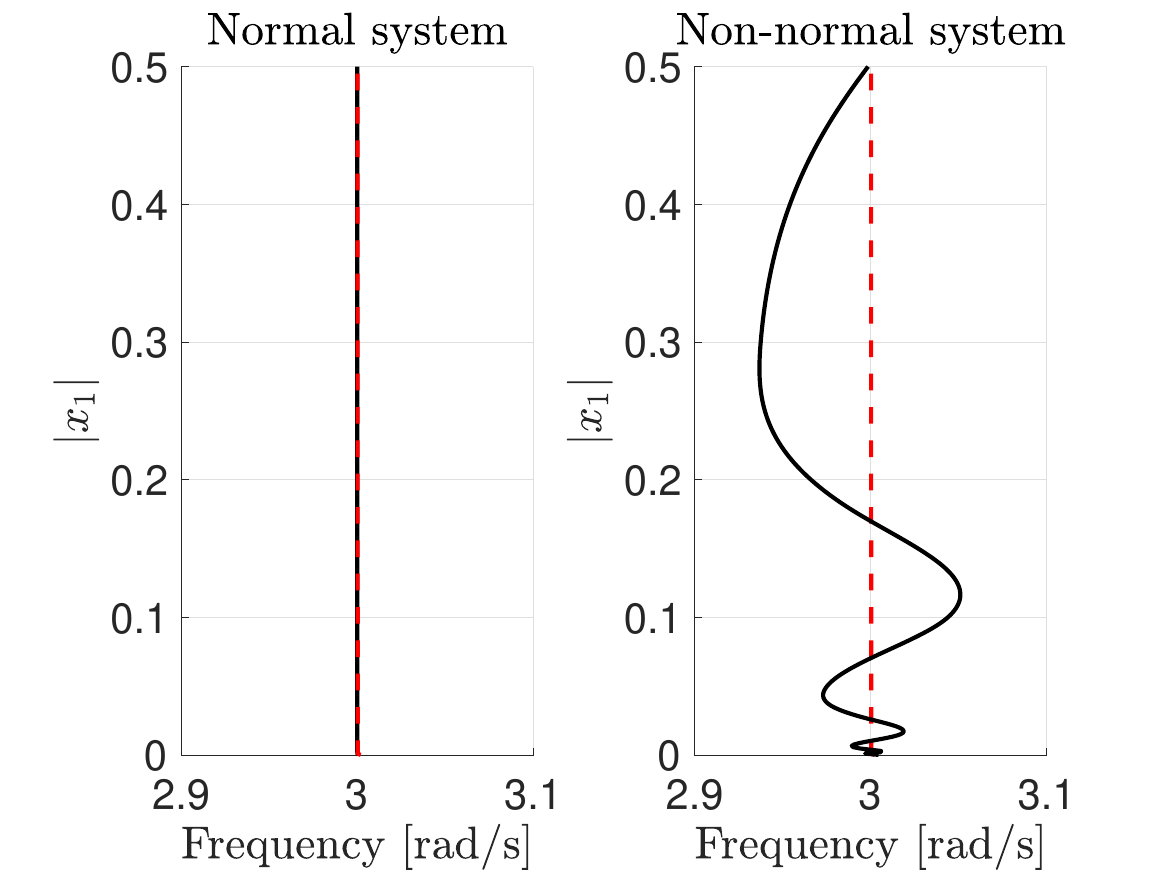}
    }
    \caption{Subfigures \subref{1d_ex_orth} and \subref{1d_ex_obl} represent the phase space of a two-dimensional linear system with real eigenvalues, when the slow and fast directions are normal (eq. \ref{2d_example}) or non-normal (eq. \ref{2d_example_nonnormal}) to each other, respectively. Subfigure \subref{2d_ex_backones}: if we consider a linear system with two oscillating modes (eq. \ref{2d_linear_system}), the backbone curve of a decaying trajectory lying on a fractional SSM (black) exhibits oscillations around the backbone curve of decaying trajectories on the primary SSM (dashed, red), when the slow and fast subspaces are not normal (non-normal linear part). Otherwise, the backbone curves are identical.} 
    \label{figures}
\end{figure*}

Suppose that we now observe the system in a specific set of coordinates $\mathbf{z}$ such that $E$ and $F$ are non-normal to each other. Namely, we introduce the coordinate change 
\begin{equation}
\mathbf{x} = \mathbf{T}\,\mathbf{z}, \quad \mathbf{T} = 
\begin{pmatrix}
1 & \frac{\delta}{\beta-\alpha} \\ 
0 & 1
\end{pmatrix},
\end{equation}
so that 
\begin{equation}\label{2d_example_nonnormal}
\dot{\mathbf{z}} = \tilde{\mathbf{A}}\mathbf{z}, \quad \tilde{\mathbf{A}} = 
\begin{pmatrix}
-\alpha & \delta \\ 0 & -\beta
\end{pmatrix}.
\end{equation}
The system is now non-normal. i.e. $\mathbf{v}^{(1),\mathrm{T}}\cdot \mathbf{v}^{(2)}\neq 0$, with 
\begin{equation}
\mathbf{v}^{(1)} = \left(1, 0 \right)^\mathrm{T},\quad \mathbf{v}^{(2)} = \left(\frac{\delta}{\alpha-\beta}, \;1 \right)^\mathrm{T}, \quad E = \text{span}\left\{ \mathbf{v}^{(1)}\right\}, \quad F = \text{span}\left\{ \mathbf{v}^{(2)}\right\}.
\end{equation}\label{2d_eigenvectors}
If we now consider an initial condition on a generic fractional SSM, its evolution in time 
\begin{equation}\label{fractional_SSM_eq_2D}
         z_1(t) = \left(z_{10} - \frac{\delta z_{20}}{\alpha - \beta} \right) e^{-\alpha t} + \frac{\delta z_{20}}{\alpha - \beta}e^{-\beta t} 
\end{equation}
cannot be modeled by the normal projection onto $\mathcal{W}\left( E\right)$, whose internal dynamics is similar to eq. \eqref{internal_SSM_2d} and reads
\begin{equation}
	z_1(t) = z_{10}e^{-\alpha t}.
\end{equation}
We observe that in eq. \eqref{fractional_SSM_eq_2D} the solution $z_1(t)$ depends on $z_{20}$, which selects a particular fractional SSM. The non-normality between eigenspaces is solely responsible for this behavior. A sketch of the phase space in the two cases is reported in Figures \ref{1d_ex_orth} and \ref{1d_ex_obl}. 

Let us consider the following system, modeling a typical damped linear oscillator, given by
\begin{equation}\label{2d_linear_system}
	\mathbf{\dot{x}} = \mathbf{A}\mathbf{x},\quad \mathbf{x} = (x_1, x_2, x_3, x_4)^\mathrm{T},\quad \mathbf{A} = 
	\begin{pmatrix}
		-\alpha & -\omega& A & B \\
		\omega & -\alpha& C & D \\ 
		0&0&-\beta & -\nu \\
		0&0&\nu& -\beta
	\end{pmatrix}, \quad\beta > \alpha > 0.
\end{equation}
If $A,B,C$ and $D$ are zero, then the two modes are orthogonal to each other (normal system), i.e. the slow spectral subspace $E$ spanned by $\left(x_1, x_2 \right)$ is orthogonal to the fast subspace $F$. Otherwise, the two modes are coupled (non-normal system). We examine decaying trajectories with initial conditions lying on the primary SSM $\mathcal{W}\left(E\right) \equiv E$ and on a fractional SSM, both with and without coupling. Given oscillatory decaying trajectories, we approximate backbone curves through the Peak Finding and Fitting (PFF) algorithm (\citet{Jin_Mengshi_2020}). For each semi-period of the signal, the PFF algorithm associates the maximum amplitude of the observed quantity ($| x_1 |$ in this example) with the corresponding local frequency, which is directly obtained from the semi-period. This process yields an amplitude-frequency pair for each semi-period of oscillation of the signal, effectively constructing points along the backbone curve. Fig. \ref{2d_ex_backones} compares the backbone curves they generate. Fractional effects, seen as oscillations in the backbone curves, appear in the non-normal system. In contrast, normal systems feature identical reduced dynamics on any SSM, whether primary or fractional. Consequently, already for a linear system, the dynamics of trajectories on fractional SSMs cannot be accurately reduced to those on the primary SSM when the system is non-normal. 

\section{Method}\label{method}
When constructing the backbone curves as above, we orthogonally project the trajectory onto the slow subspace and then apply the PFF procedure in order to identify the history of instantaneous frequencies and amplitudes along decaying oscillations.

\begin{figure*}[b]
    \centering
    \subfloat[\label{normal_vs_oblique_backbone}]{
        \includegraphics[scale=0.26]{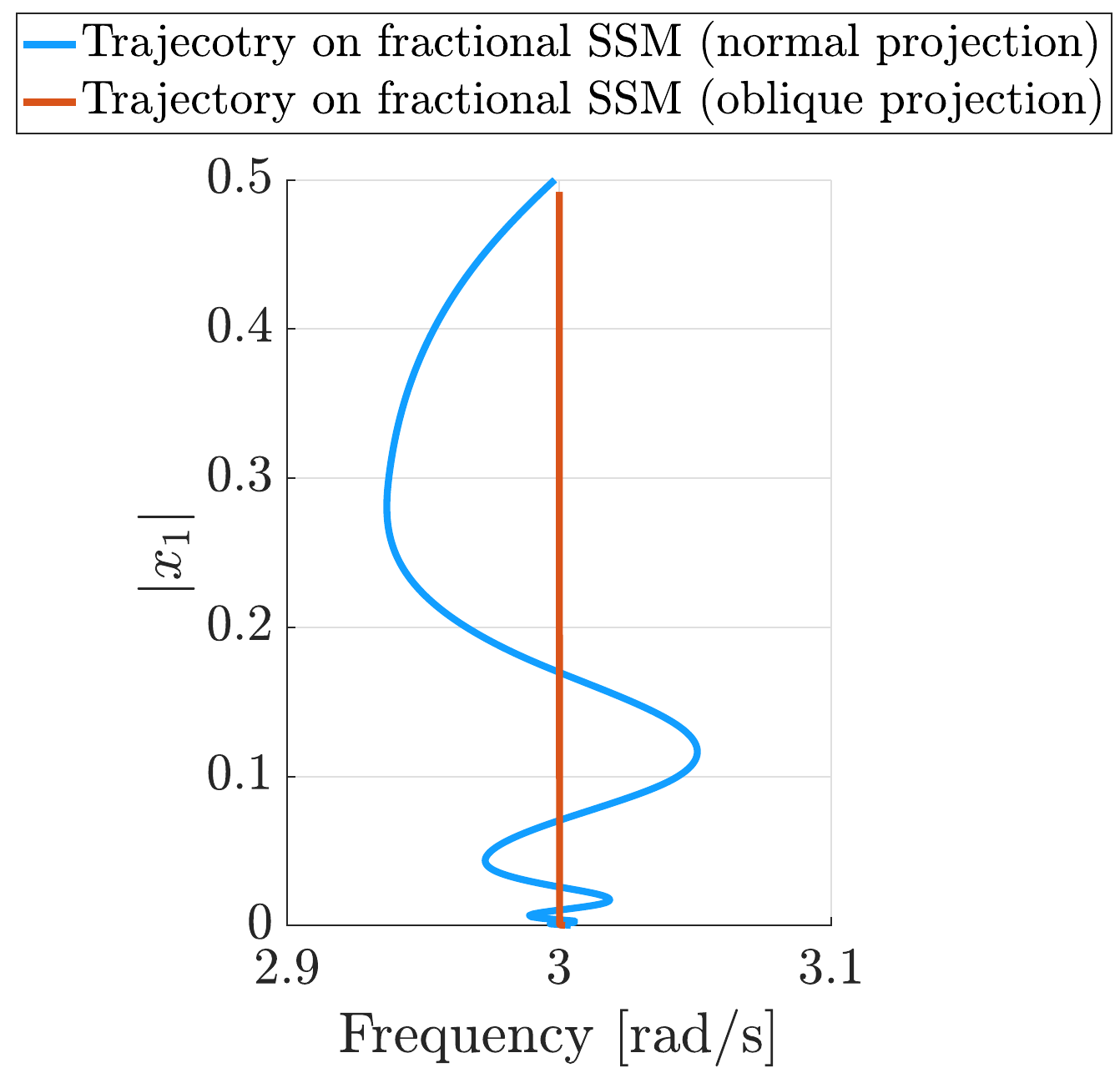}
    }
    \quad
    \subfloat[\label{figure_oblique_projection}]{
        \includegraphics[scale=0.36]{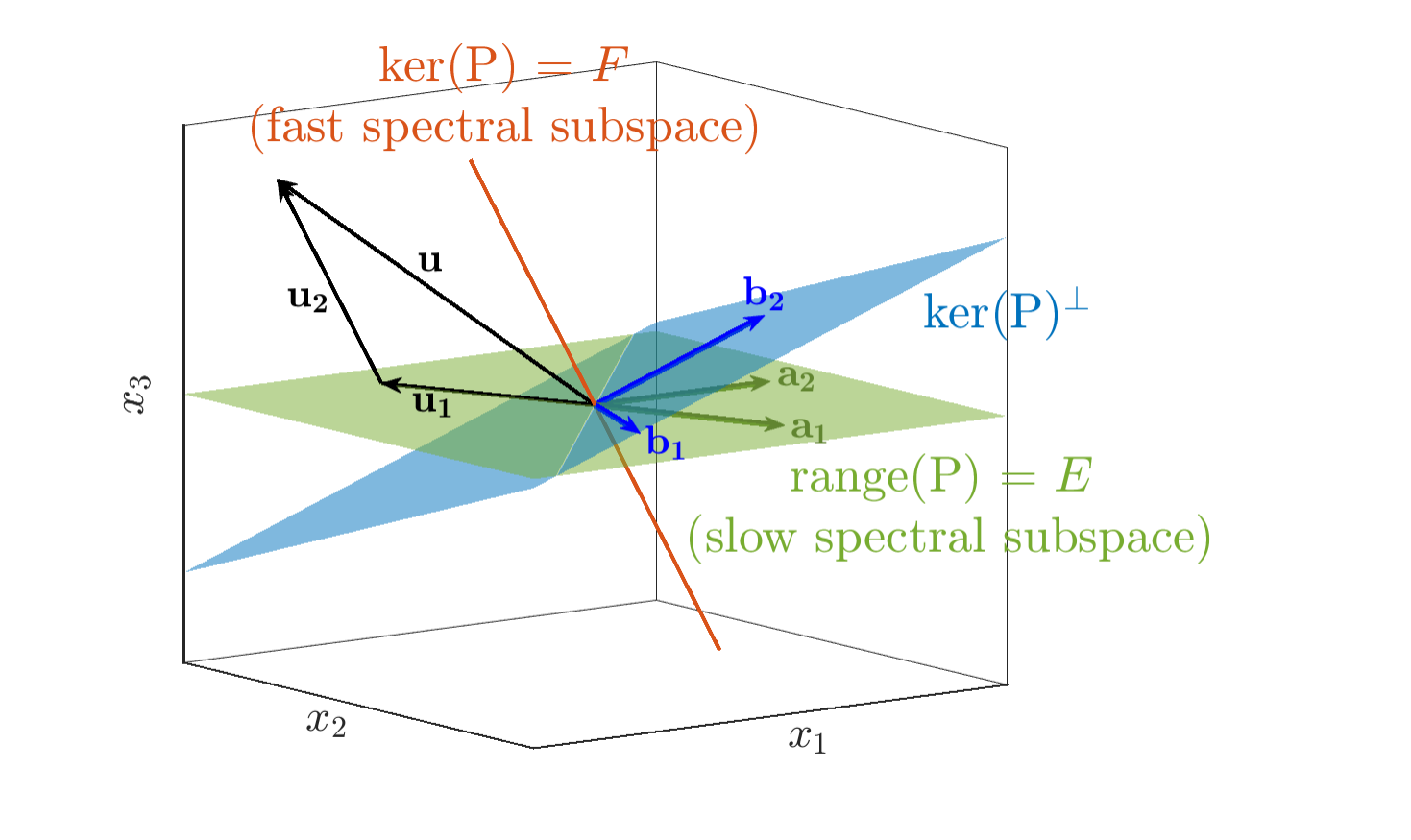}
    }
    \caption{Aspects of non-normality between the slow spectral subspace $E$ and the fast spectral subspace. \subref{normal_vs_oblique_backbone} Two different ways of extracting the backbone curve of system \eqref{2d_linear_system} with parameter values $\alpha = 0.3$, $\beta = 0.63$, $\omega = 3$, $\nu = 8$, $A = B = C = D = 1$, $\mathbf{x}_0 = \left(1,1,0.8,0.8 \right)^\mathrm{T}$: via normal projection onto $E$ (light blue) and via oblique, $F$-parallel projection onto $E$ (red). \subref{figure_oblique_projection} The construction of an $F$-parallel oblique projection $\mathbf{P}$ onto $E$.} 
\label{fig_normal_vs_oblique_projection}
\end{figure*}

We now discuss what happens if we project onto the slow subspace in an oblique way, along the direction of the fast subspace, rather than normally. Notably, this results in a straight backbone curve, which is identical to the one of the primary SSM (see Fig. \ref{normal_vs_oblique_backbone}). Consequently, the reduced dynamics on the primary SSM can serve now as a reduced-order model for the dynamics on fractional SSMs, as long as trajectory positions are projected obliquely onto $E$ along the direction of $F$. 

We now recall the general construction of an oblique projection $\mathbf{P}: \mathbb{R}^n \to E$ for system \eqref{eq:basic_eq} (see \citet{Banerjee_2014} for a full explanation). Referring to Fig. \ref{figure_oblique_projection}, we decompose a general vector $\mathbf{u} \in \mathbb{R}^n$
\begin{equation}\label{general_vector}
\mathbf{u} = \mathbf{u_1} + \mathbf{u_2}, \quad  \mathbf{u_1}\in E, \quad  \mathbf{u_2}\in F.
\end{equation}
For the envisioned projection $\mathbf{P}$, we have 
\begin{equation}\label{range_ker_projection}
	E = \text{range}\left( \mathbf{P}\right), \quad F = \text{ker}\left( \mathbf{P}\right).
\end{equation}
By equations \eqref{general_vector}-\eqref{range_ker_projection}, $\mathbf{u_1}$ is in the range of the projection and hence 
\begin{equation}\label{u_1_w}
\mathbf{u_1}  = \mathbf{Q}\mathbf{w},\quad  \mathbf{Q} \in \mathbb{R}^{n\times d}, \quad \mathbf{w} \in \mathbb{R}^d, \quad \text{dim}\left(E\right) = d.
\end{equation}
Here, the columns of $\mathbf{Q} \in \mathbb{R}^{n \times d}$ form a basis in $E$. For instance, for $d = 2$, we have $\mathbf{Q} = \left[\mathbf{a}_1 | \mathbf{a}_2 \right]$, with the vectors $\mathbf{a}_1, \mathbf{a}_2 \in E$ shown in Fig. \ref{figure_oblique_projection}. Furthermore, we have 
\begin{equation}\label{B_equation}
\mathbf{B}^\mathrm{T}\mathbf{u_2} = \mathbf{B}^\mathrm{T}\left(\mathbf{u}-\mathbf{u_1} \right) = \mathbf{0},
\end{equation}
where $\mathbf{B} \in \mathbb{R}^{n\times d}$ is a matrix whose columns form a basis of $F^\perp$, the orthogonal complement of $F = \text{ker} \left(\mathbf{P}\right)$.  For instance, for $d = 2$, we have $\mathbf{B} = \left[\mathbf{b}_1 | \mathbf{b}_2 \right]$, with the vectors $\mathbf{b}_1, \mathbf{b}_2 \in F^\perp$ shown in Fig. \ref{figure_oblique_projection}.

Substituting eq. \eqref{u_1_w} into \eqref{B_equation}, we obtain
\begin{equation}\label{formula_w}
\mathbf{w} = \left(\mathbf{B}^\mathrm{T}\mathbf{Q} \right)^{-1}\mathbf{B}^\mathrm{T} \mathbf{u}.
\end{equation}
Substituting formula \eqref{formula_w} into \eqref{u_1_w} then gives $\mathbf{u_1} = \mathbf{P}\mathbf{u}$ with 
\begin{equation}\label{oblique_projection}
	\mathbf{P} = \mathbf{Q}\left(\mathbf{B}^\mathrm{T}\mathbf{Q} \right)^{-1}\mathbf{B}^\mathrm{T}.
\end{equation}
In the special case of $E \perp F$, $\mathbf{P}$ is a normal projection, its range coincides with the orthogonal complement of its kernel, namely $\mathbf{Q} \equiv \mathbf{B}$ and 
\begin{equation}
	\mathbf{P} = \mathbf{P}^\mathrm{T}=  \mathbf{Q}\left(\mathbf{Q}^\mathrm{T}\mathbf{Q} \right)^{-1}\mathbf{Q}^\mathrm{T}.
\end{equation}

Formula \eqref{oblique_projection} gives the general form of oblique projections onto $E$. SSM-reduction algorithms approximated $E$ as the span of the dominant singular vectors of the trajectory data matrix. For our purposes here, we need a more accurate approximation, which prompts us to use a combination of singular value decomposition (SVD) with the dynamic mode decomposition (DMD).

It remains to determine the matrix $\mathbf{B}$ in formula \eqref{oblique_projection} to obtain the projection matrix $\mathbf{P}$. As the columns of $\mathbf{B}$ span $F^\perp$, identifying $\mathbf{B}$ from data is a much easier problem than identifying the (generally very high dimensional) subspace $F$ spanned by the fast eigenspaces. We find $\mathbf{B}$ as the matrix that renders the variance of the backbone curve minimal when computed from the input data after the application of the projection $\mathbf{P}$. The steps in this procedure are as follows.

\begin{enumerate}
\item \textbf{Linear regime identification.}\label{step_1} We identify trajectory pieces falling in the linear regime of the system by looking at the frequency values of the backbone curve and computing the average as we remove one frequency at a time. The amplitude range corresponding to the frequencies whose difference in average falls below a user-defined threshold identifies the linear regime.

\item \textbf{2D slow subspace identification.}\label{step_2} We aim to identify the matrix $\mathbf{Q} = \left[\mathbf{a}_1 | \mathbf{a}_2 \right]$ spanning the 2D slow subspace $E$ which is required to construct the oblique projection \eqref{oblique_projection}. The subspace $E$ is tangent to the primary SSM that we want to reconstruct (\citet{haller2016}) and it contains the reduced coordinates that parametrize it. Practically, we employ DMD on the linear regime of the data (\citet{Schmid_2008}) along with SVD in order to approximate $E$. The input data to this procedure consists of N snapshots of a subset of the system's states collected in the observable vector $\mathbf{y} \in \mathbb{R}^p$ with $p\leq n$, collected in the matrices
\begin{equation}
\begin{aligned}
&
\mathbf{V}_1 = 
\begin{pmatrix}
\mathbf{y}(t_1), & \mathbf{y}(t_2) ,& \dots, & \mathbf{y}(t_N)
\end{pmatrix} \quad \in \mathbb{R}^{p\times N}
\\ &
\mathbf{V}_2 = 
\begin{pmatrix}
\mathbf{y}(t_2), & \mathbf{y}(t_3) ,& \dots, & \mathbf{y}(t_{N+1})
\end{pmatrix} \quad \in \mathbb{R}^{p\times N},
\end{aligned}
\end{equation}
so that 
\begin{equation}\label{equation_DMD}
\mathbf{V}_2 = \mathbf{A} \mathbf{V}_1,
\end{equation}
where $\mathbf{A} \in \mathbb{R}^{p\times p}$ is an appropriate matrix. Applying the truncated compact version of SVD (\citet{Banerjee_2014}) to the snapshot matrix $\mathbf{V}_1$ yields 
\begin{equation}\label{svd_equation}
\mathbf{V}_1 = \mathbf{U} \boldsymbol{\Sigma}\mathbf{W}^\mathrm{T},
\end{equation}
where $\mathbf{U} \in \mathbb{R}^{p \times d}$, $\boldsymbol{\Sigma} \in \mathbb{R}^{d \times d}$ and $\mathbf{W} \in \mathbb{R}^{N \times d}$, so that $\text{det}\left(\boldsymbol{\Sigma} \right) \neq 0$ and $\mathbf{U}^\mathrm{T}\mathbf{U} = \mathbf{W}^\mathrm{T}\mathbf{W} = \mathbf{I}_d$. 

We now substitute eq. \eqref{svd_equation} into eq. \eqref{equation_DMD} and premultiply by $\mathbf{U}^\mathrm{T}$ to obtain
\begin{equation}
	 \mathbf{U}^\mathrm{T}\mathbf{V}_2 = \mathbf{U}^\mathrm{T}\mathrm{A}\mathbf{U}\boldsymbol{\Sigma}\mathbf{W}^\mathrm{T}
\end{equation}
and hence 
\begin{equation}
\mathbf{U}^\mathrm{T}\mathbf{A}\mathbf{U} = \mathbf{U}^\mathrm{T}\mathbf{V}_2\mathbf{W}\boldsymbol{\Sigma}^{-1} = \mathbf{\tilde{A}}, \qquad \mathbf{\tilde{A}} \in \mathbb{R}^{d \times d}.
\end{equation}
The eigenvectors of $\mathbf{\tilde{A}}$ are readily computable and collected as columns of matrix $\mathbf{V}_{\tilde{A}} \in \mathbb{R}^{d\times d}$. The matrix $\mathbf{Q} \in \mathbb{R}^{p \times d}$ that approximates the slow subspace $E$ is found as the matrix of eigenvectors of $\mathbf{A}$, i.e. $\mathbf{Q} = \mathbf{U} \mathbf{V}_{\tilde{A}}$.

\item \textbf{Matrix B identification.}\label{step_3} We find $\mathbf{B}$ as the solution of the optimization problem 
\begin{equation}
\label{eq:constfunction}
\mathbf{B}^*= \underset{\mathbf{B} \in \mathbb{R}^{p \times 2}}{\arg\min} \,\sigma \left(\mathbf{Q}\left(\mathbf{B}^\mathrm{T}\mathbf{Q}\right)^{-1}\mathbf{B}^\mathrm{T}\,\boldsymbol{y}_\mathrm{lin} \right),
\end{equation}
where $\boldsymbol{y}_\mathrm{lin}$ denotes the data in the linear regime and $\sigma \left( \cdot \right)$ measures the variance of the backbone curve of its argument. The optimization problem is addressed using a Quasi-Newton algorithm, where $\mathbf{Q}$ serves as the initial condition for convergence to $\mathbf{B}$.

Once we have determined the matrices $\mathbf{Q}$ and $\mathbf{B}$ from steps \ref{step_1}-\ref{step_3}, we compute the oblique projection $\mathbf{P}$ from eq. \eqref{oblique_projection}. This defines the projected coordinates on $E$ as
\begin{equation}\label{red_coord_projection}
\mathbf{z} = \mathbf{P}\mathbf{y},
\end{equation}
for each column $\mathbf{y}$ of the data matrices $\mathbf{V}_1$ and $\mathbf{V}_2$ that lie close to the primary SSM.

\item \textbf{SSM parametrization and reduced dynamics.}\label{step_4} We now follow a procedure similar to that of \citet{cenedese2022}. Specifically, we approximate the parametrization of the SSM and its reduced dynamics with multivariate polynomials of the reduced coordinates $\boldsymbol{\xi}$, defined as 
\begin{equation}
\boldsymbol{\xi} = \mathbf{\tilde{Q}}^\mathrm{T}\mathbf{z},
\end{equation}
where matrix $\mathbf{\tilde{Q}}$ is an orthonormal representation of the tangent space $E$.
The parametrization is sought in the form of a multivariate polynomial
\begin{equation}
\mathbf{y} = \mathbf{h}\left(\boldsymbol{\xi} \right) = \mathbf{M}_1 \boldsymbol{\xi}  + \mathbf{M}\boldsymbol{\xi}^{2:M}
\end{equation}

The unknown $\mathbf{M}_1$ and $\mathbf{M}$ are found via constrained regression by the following minimization problem:

\begin{equation}
\left(\mathbf{M}_1^*, \mathbf{M}^*\right) = \underset{\mathbf{M}_1, \mathbf{M}}{\arg\min} \left \| \mathbf{y} - \mathbf{M}_1\boldsymbol{\xi} - \mathbf{M}\left(\boldsymbol{\xi} \right)^{2:M} \right\|^2.
\end{equation}
\end{enumerate}

\begin{figure}[h!]
    \centering
    \includegraphics[width = 0.5\textwidth]{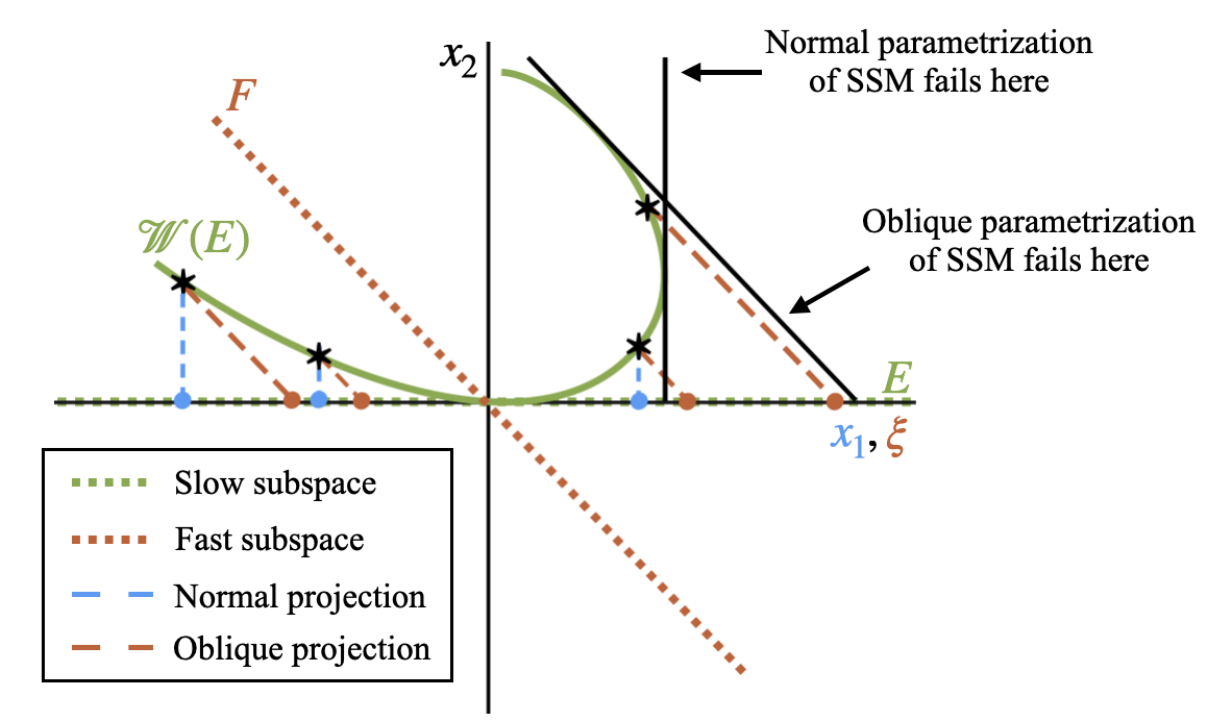}
    \caption{Comparison between normally and obliquely projected reduced coordinates that parametrize the SSM. We note that the $\xi$ coordinate obtained from the oblique projection parametrizes $\mathcal{W}\left(E\right)$ on a much larger domain than the $x_1$ coordinate obtained from normal projection.}
    \label{ssm_parametrization}
\end{figure}

To clarify the advantage of the oblique projection as constructed in Step \ref{step_3} for the SSM parametrization, let us consider a general 2D nonlinear system, with the observables $\mathbf{y} = \left(x_1, x_2 \right)^\mathrm{T}$, where the coordinate $x_1$ lies in the slow spectral subspace $E$ of the linearized system. We seek a parametrization of the SSM $\mathcal{W}\left(E\right)$ as a function of the reduced coordinate $\xi$ in the form 
\begin{equation}\label{SSM_param_red}
	x_2 = h\left(\xi \right).
\end{equation}
If $\xi \equiv x_1$, then the reduced coordinate is the normal projection onto the slow subspace (blue, dashed lines in Fig. \ref{ssm_parametrization}). In contrast, the oblique projection $\mathbf{P}$ is defined by the relationship $\xi = p\left(x_1, x_2\right)$, which gives an SSM parametrization in the implicit form
\begin{equation}\label{SSM_param_implicit}
	x_2 = h \Bigl( p \left( x_1, x_2 \right) \Bigr),
\end{equation}
visualized as the red, dashed lines in Fig. \ref{ssm_parametrization}. We note that once we fix the value of the reduced coordinate $\xi$, also $x_2$ on the SSM is fixed by the relationship \eqref{SSM_param_red}, while we have to solve eq. \eqref{SSM_param_implicit} in order to retrieve the corresponding $x_1$.

Steps \ref{step_1}-\ref{step_4} have been implemented in a publicly available, updated version of the SSMLearn algorithm, which is downloadable from \url{https://github.com/haller-group/SSMLearn}.

\begin{figure}[H]
    \centering
    \subfloat[\label{2d_example_SSMLearn_orthogonal}]{
        \includegraphics[scale=0.33]{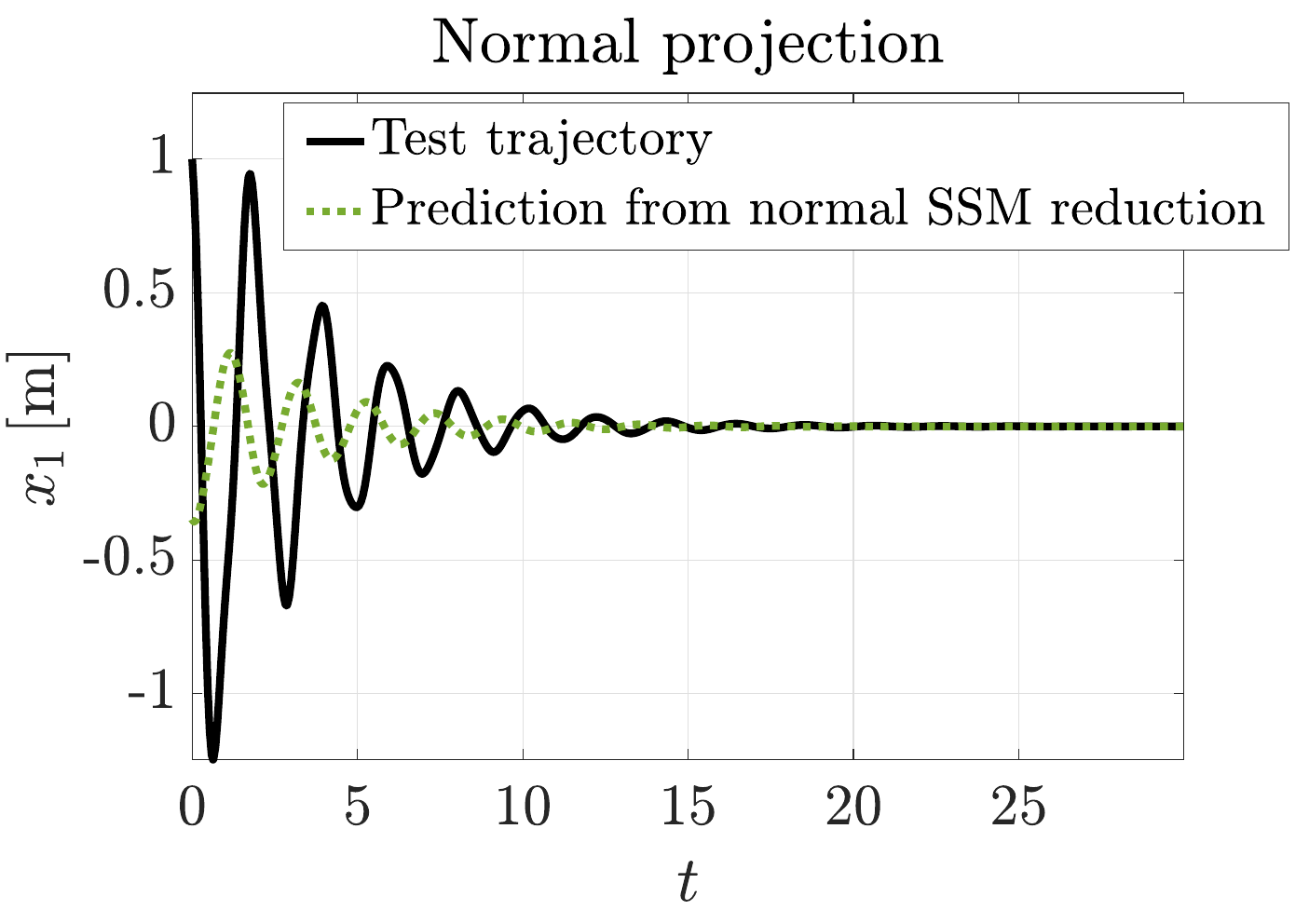}
    }
    \quad
    \subfloat[\label{2d_example_SSMLearn_oblique}]{
        \includegraphics[scale=0.33]{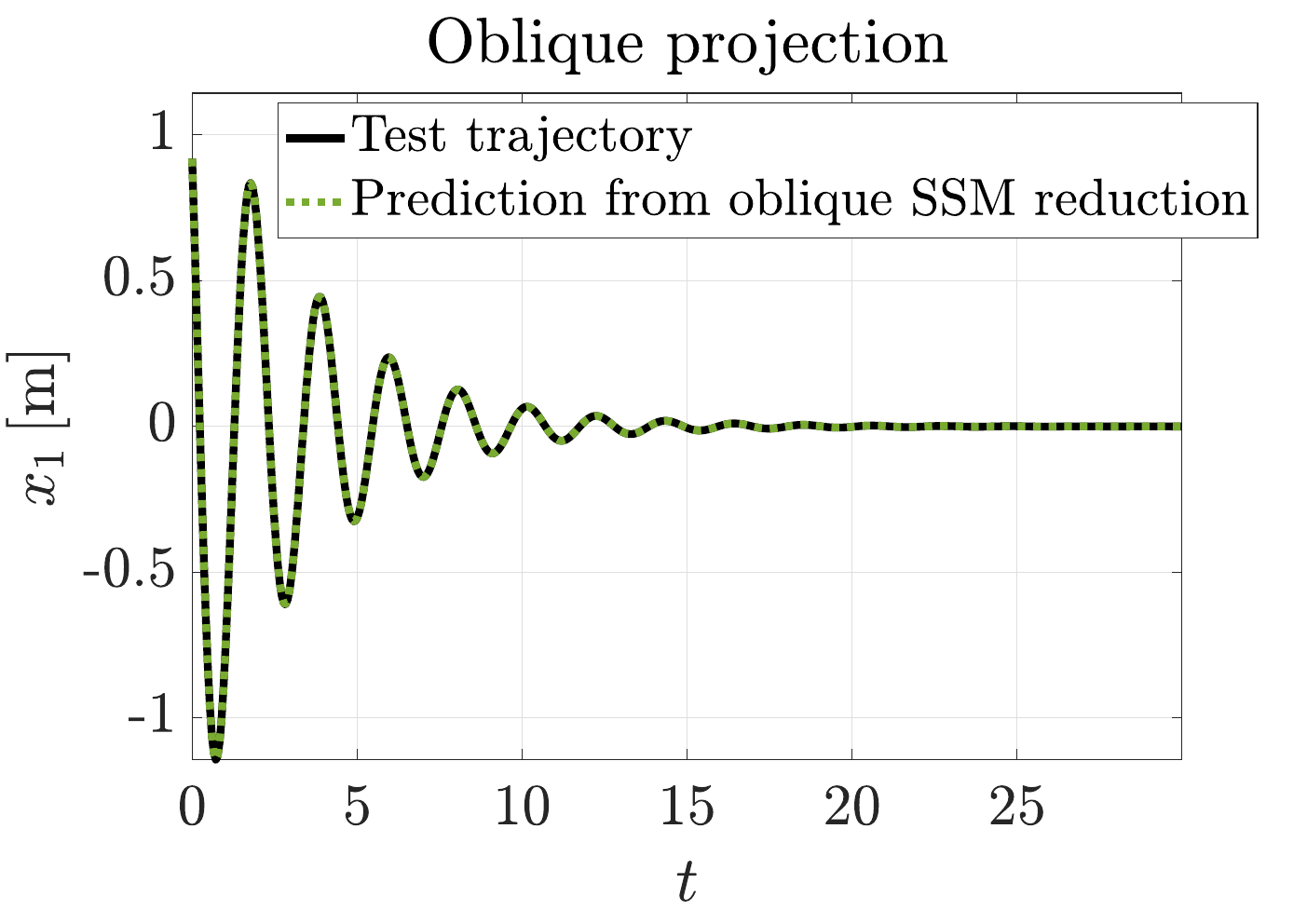}
    }
    \caption{Decaying test trajectory of the linear non-normal system \eqref{2d_linear_system}, with initial condition on a fractional SSM and its prediction from \subref{2d_example_SSMLearn_orthogonal} normal SSM reduction and \subref{2d_example_SSMLearn_oblique} from oblique SSM reduction. Both SSMs are computed from data using the SSMLearn algorithm \url{https://github.com/haller-group/SSMLearn}.}  
    \label{2d_example_SSMLearn}
\end{figure}

As a simple illustration of the above procedure, we now apply oblique SSM-reduction to the linear non-normal system \eqref{2d_linear_system}. As seen in Fig. \ref{2d_example_SSMLearn}, the internal dynamics of the primary SSM $\mathcal{W}\left(E\right)$ parametrized by the reduced coordinate $\xi$ successfully tracks the decaying trajectory which gave rise to the oscillatory backbone curve in Fig. \ref{2d_ex_backones}, as opposed to the standard normal-constructed primary SSM (\ref{2d_example_SSMLearn_orthogonal}). 
Both cases employ the SSMLearn algorithm developed in \citet{cenedese2022} for computing the SSMs through a data-driven approach, equipped with oblique projection when it is used. The oblique projection is determined as the one that minimizes the oscillations of the backbone curve. 

\section{Results}\label{results}
Next, we apply the oblique SSM reduction to several numerical and experimental examples. We show that this method proves to be effective when a simple normal SSM reduction fails.

\subsection{Shaw-Pierre oscillator chain}\label{SP_example}
\begin{figure*}[b]
	\centering
	\includegraphics[scale=1.3]{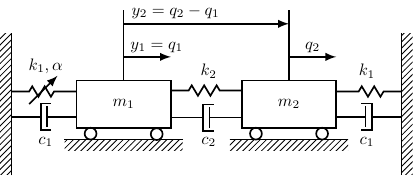}
	\caption{Geometry of the two-degree-of-freedom mechanical system studied in \citet{shaw_pierre_93}, modified by adding a damper between the left mass and the wall.}
	\label{SP_system}
\end{figure*}
We consider the two-degree-of-freedom mechanical system studied by \citet{shaw_pierre_93}, with an additional damper connecting the left mass to the wall. Shown in Fig. \ref{SP_system}, the system was originally studied in the $\left(q_1, q_2 \right)$ coordinates, whereas here we use the coordinates
\begin{equation}
	y_1 = q_1, \quad y_2 = q_2 - q_1.
\end{equation}
The governing equations are
\begin{equation}\label{mechanical_system_equation}
\mathbf{M}\mathbf{\ddot{y}} + \mathbf{C}\mathbf{\dot{y}} + \mathbf{K}\mathbf{y} + \mathbf{f}_\mathrm{nl} = \mathbf{0},
\end{equation}
with $ \mathbf{y} = \left(y_1,\, y_2 \right)^\mathrm{T}$ and 
\begin{equation}\label{SP_governing_equations}
\begin{aligned}
\mathbf{M} = 
   \begin{pmatrix}
    m_1 & 0\\
    m_2 & m_2
    \end{pmatrix},\quad
 \mathbf{C} = 
   \begin{pmatrix}
    c_1 & -c_2\\
    c_1 & c_1 + c_2
    \end{pmatrix},\quad
\mathbf{K} = 
   \begin{pmatrix}
    k_1 & -k_2\\
    k _1& k_1 + k_2
    \end{pmatrix},\quad
\mathbf{f}_\mathrm{nl} = 
    \begin{pmatrix}
    -\alpha y_1^3 \\
    0
    \end{pmatrix}.
\end{aligned}
\end{equation}

\begin{figure*}[]
    \centering
    \subfloat[\label{SP_manifold_with_trajectories}]{
        \includegraphics[scale=0.27]{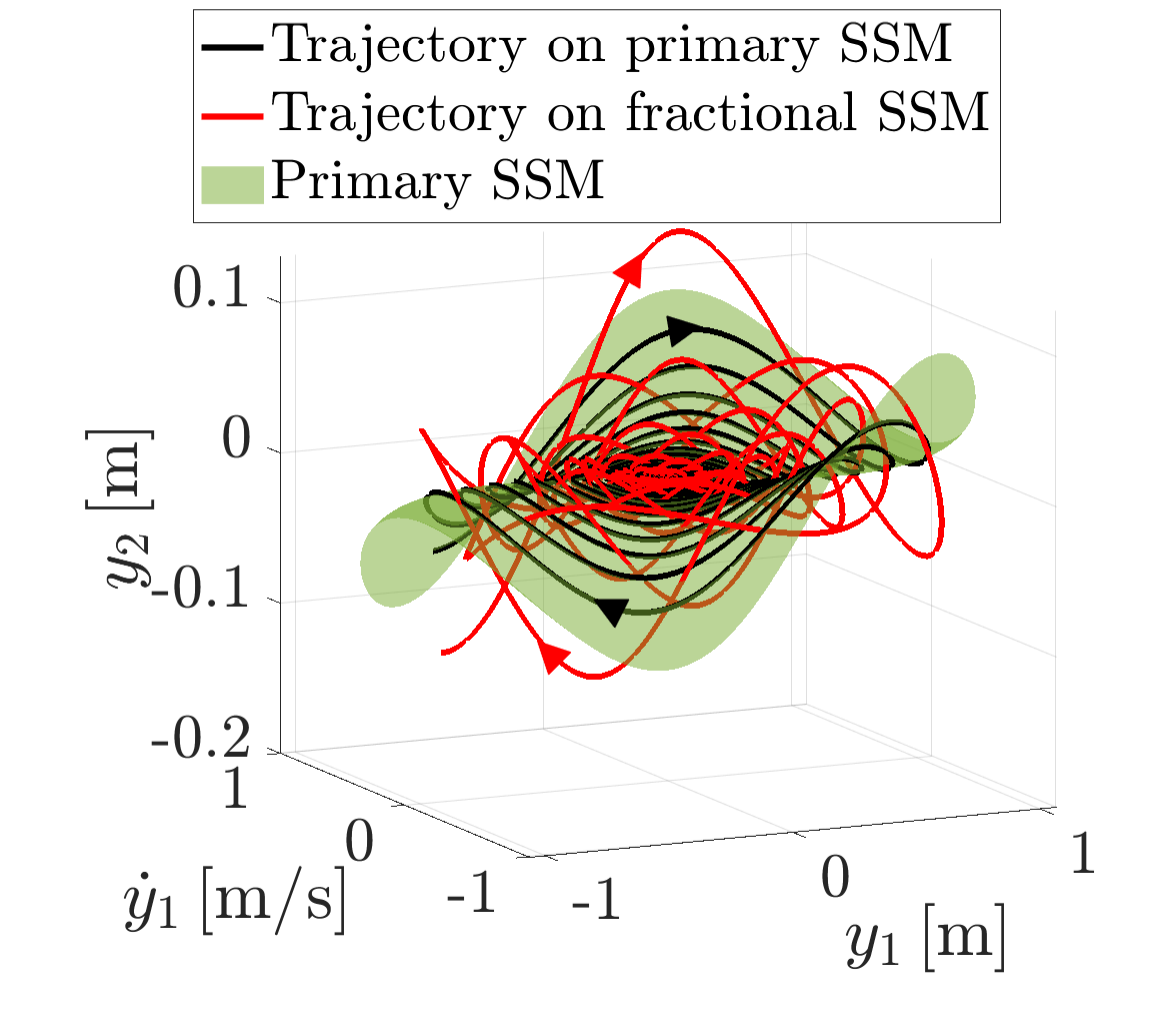}
    }\quad
    \subfloat[\label{SP_backbone_curves}]{
        \includegraphics[scale=0.18]{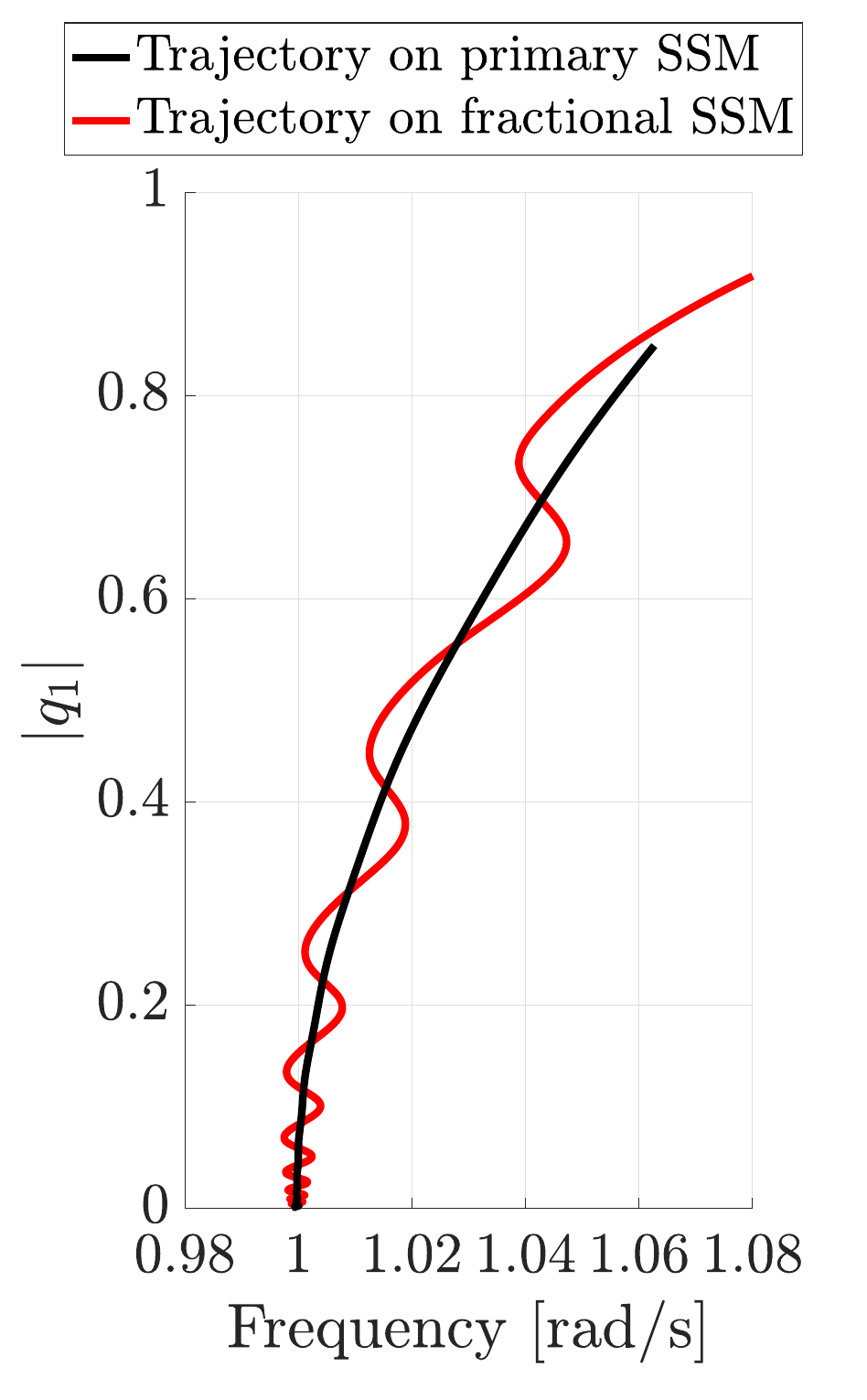}
    }
    \quad
    \subfloat[\label{SP_linear_regime}]{
        \includegraphics[scale=0.18]{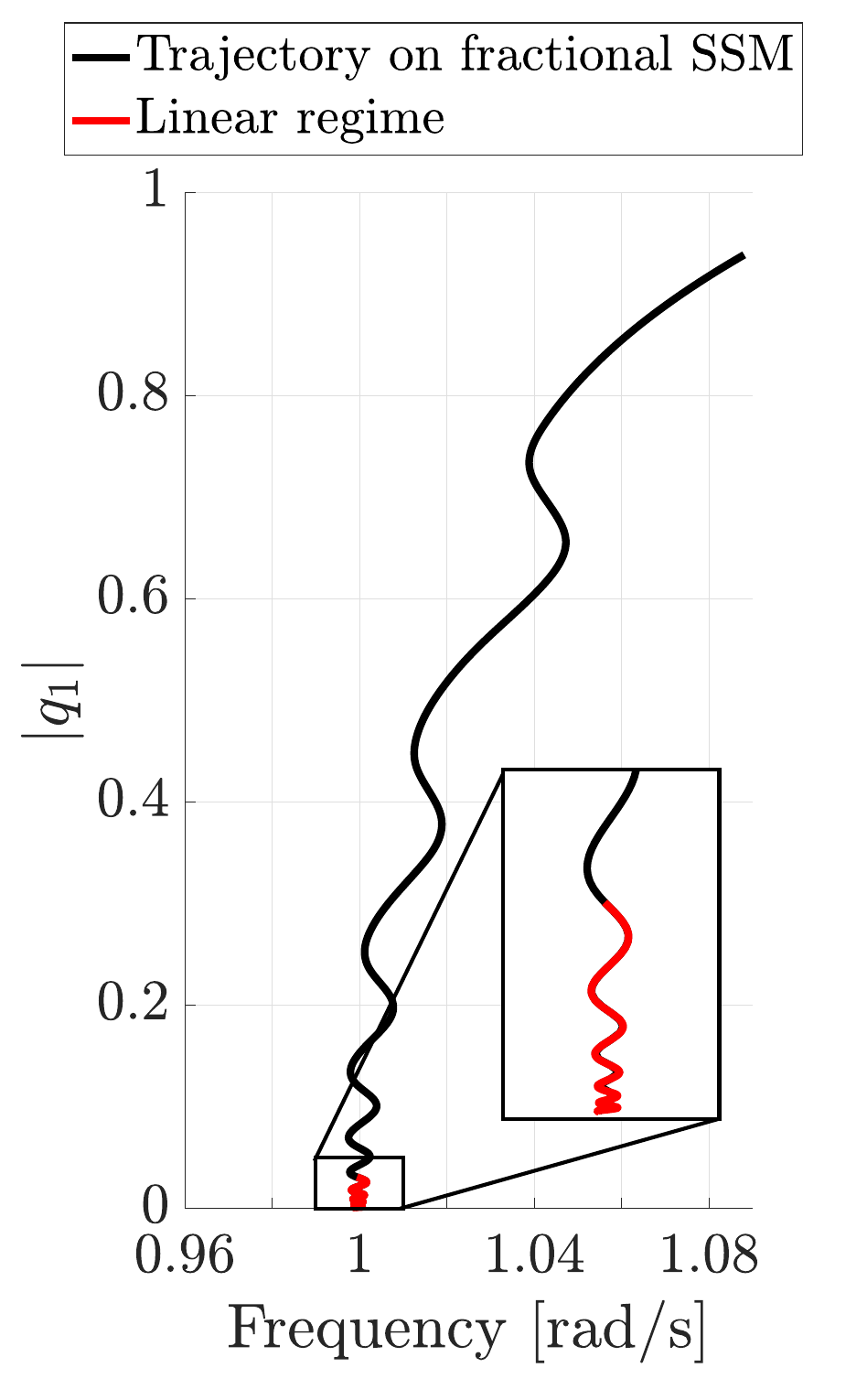}
    } 
    \caption{Trajectories of the Shaw-Pierre system observed in the $\left(y_1, y_2 \right)$ coordinates. \subref{SP_manifold_with_trajectories} The primary SSM $\mathcal{W}\left(E\right)$ (green), a trajectory on it (black) and a trajectory outside $\mathcal{W}\left(E\right)$ on a fractional SSM (red). \subref{SP_backbone_curves} Backbone curve computed from observations of these trajectories. \subref{SP_linear_regime} The linear data regime to be used to construct the oblique projection $\mathbf{P}$.} 
    \label{SP_primary_fractional_trajectories}
\end{figure*}

With the notation $x_1 = y_1$, $x_2 = y_2$, $x_3 = \dot{y}_1$ and $x_4 = \dot{y}_2$, we can rewrite eq.\eqref{mechanical_system_equation} as a first-order system of ODEs
\begin{equation}\label{eq:x_coordinates}
    \dot {\mathbf{x}} =
    \mathbf{A}\mathbf{x} + \mathbf{F}_\mathrm{nl},
\end{equation}
with
\begin{align*}
    &
    \mathbf{A} = \begin{pmatrix}
    0 & 0 & 1 & 0 \\
    0 & 0 & 0 & 1 \\ \displaystyle
    -\frac{k_1}{m_1} &\displaystyle \frac{k_2}{m_1} &\displaystyle -\frac{c_1}{m_1} &\displaystyle  \frac{c_2}{m_1} \\ \displaystyle
    -\frac{\left(m_1 - m_2 \right) k_1}{m_1 m_2} &\displaystyle -\frac{m_1 \left( k_1 + k_2 \right) + m_2 k_2}{m_1 m_2} &\displaystyle -\frac{\left(m_1 - m_2 \right) c_1}{m_1 m_2} &\displaystyle -\frac{m_1 \left( c_1 + c_2 \right) + m_2 c_2}{m_1 m_2}
       \end{pmatrix},
    \\[\medskipamount] &
    \mathbf{F}_\mathrm{nl} = \begin{pmatrix}
    0 \\
    0\\\displaystyle -\frac{\alpha}{m_1} x_1^3 \\0
    \end{pmatrix}.
\end{align*}
As parameter values, we take 
\begin{equation*}
    m_1 = m_2 = 1, \quad c_1 = 0.05, \quad c_2 = 0.01, \quad k_1 = 1, \quad k_2 = 3.325, \quad \alpha = 0.5.
\end{equation*}
The eigenvalues of $\mathbf{A}$ are 
\begin{equation*}
    \begin{array}{l}
          \lambda_{1,2} = -0.025 \pm i\, 0.9997,\\
          \lambda_{3,4} = -0.035\pm i\,2.7656 ,
    \end{array}
\end{equation*}
whose eigenvectors give rise to two 2D real invariant subspaces, the slow spectral subspace $E$ and the fast spectral subspace $F$. According to our choice of coordinates, $E$ and $F$ are not normal to each other.
A decaying trajectory with initial condition outside the primary SSM $\mathcal{W}\left(E\right)$ (shown in red in Fig. \ref{SP_manifold_with_trajectories}), generates an oscillatory backbone curve shown in red in Fig. \ref{SP_backbone_curves}. This phenomenon occurs even at low amplitudes, when the dynamics are approximately linear. This is the effect of non-normality on trajectories lying on fractional SSMs outside $\mathcal{W}\left(E\right)$, as we discussed in earlier sections.

We now perform an oblique SSM reduction to the 2D slow SSM $\mathcal{W}\left(E\right)$ using the procedure outlined in section \ref{method}.
With reference to the notation used in this section, we have the dimension of the phase space $n = 4$, the dimension of the reduced-order model $d = 2$ and that of the observable space $p = 8$. We compute both normal and oblique SSMs up to $5^{th}$ order using the SSMLearn algorithm of \citet{cenedese2022} from a single trajectory, first with normal projection onto $E$ (normal SSM reduction), then with oblique projection (oblique SSM reduction).
The reduced-order model in polar coordinates reads
\begin{equation}
\begin{cases}
\dot{\rho} = -0.025 \rho + 0.017 \rho^3 - 0.091 \rho^5 \\
\dot{\theta} = 1 + 0.245 \rho^2 - 0.163 \rho^4
\end{cases}.
\end{equation}

\begin{figure*}[b]
    \centering
    \subfloat[\label{SP_orthogonal_FRC}]{
        \includegraphics[scale=0.3]{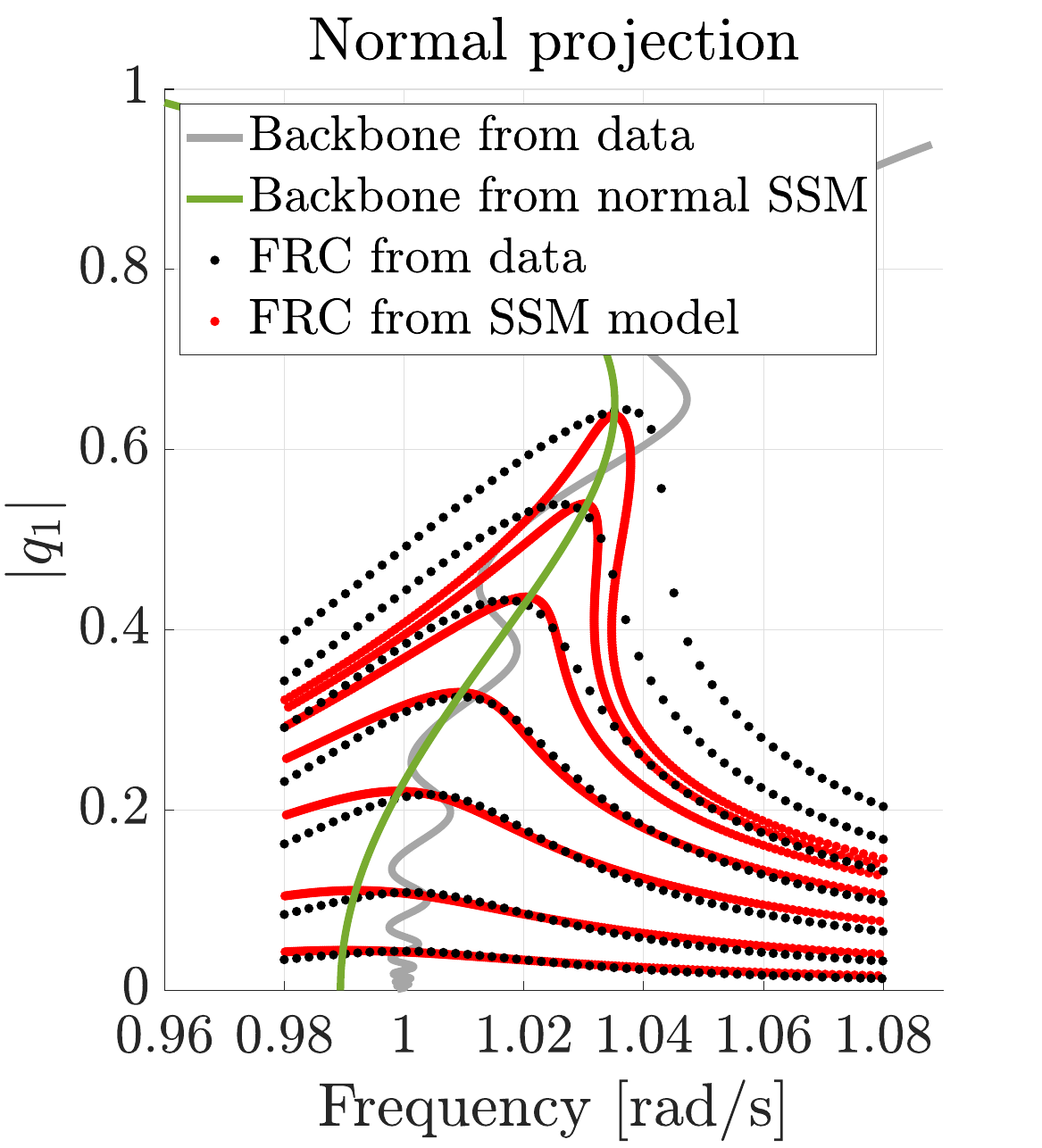}
    }
    \quad
    \subfloat[\label{SP_oblique_FRC}]{
        \includegraphics[scale=0.3]{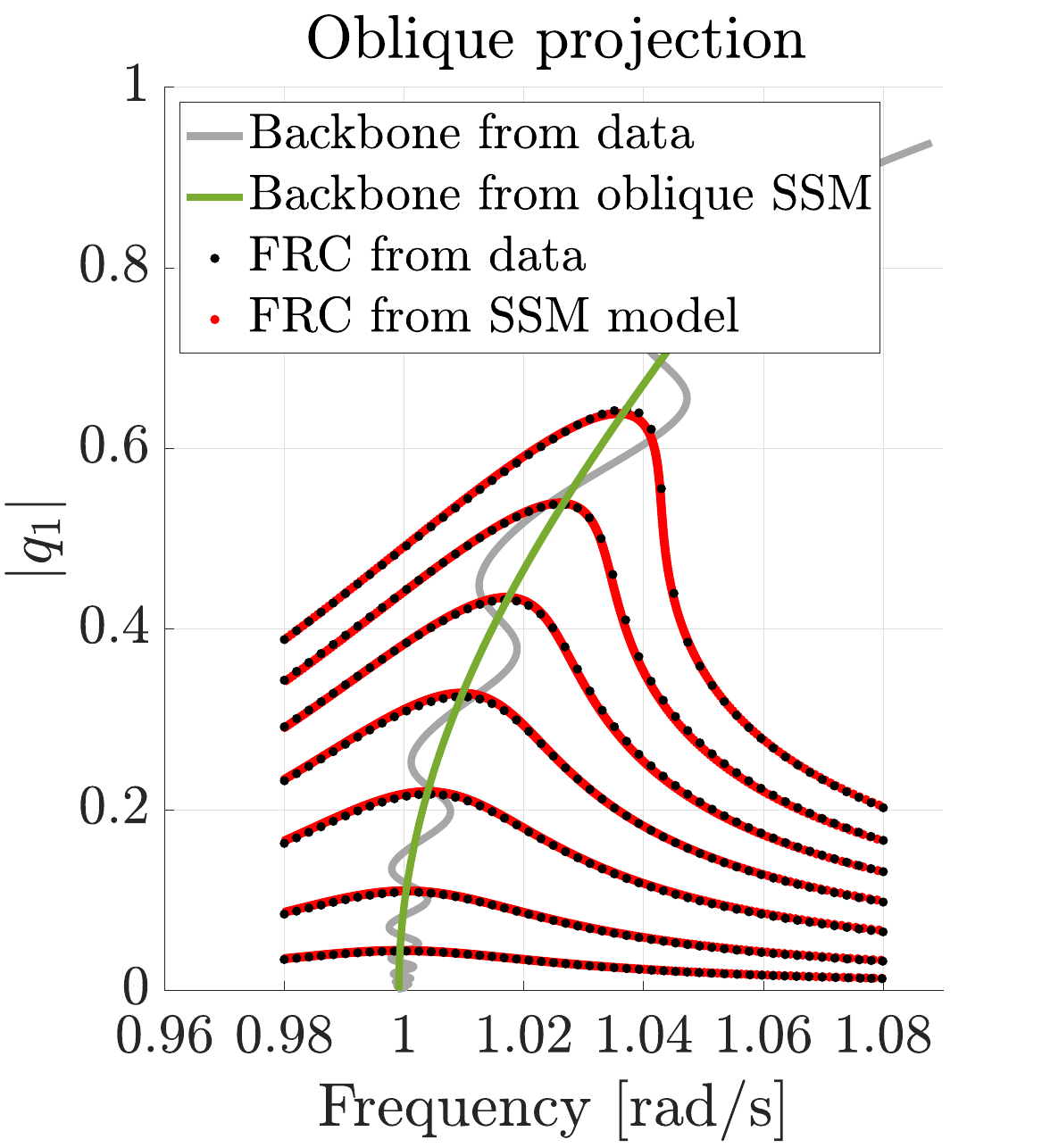}
    } 
    \caption{Predictions for backbone curve and forced response curves (FRC) in system \eqref{mechanical_system_equation}. \subref{SP_orthogonal_FRC} Normal SSM reduction. \subref{SP_oblique_FRC} Oblique SSM reduction.} 
    \label{SP_results}
\end{figure*}

A comparison between the predictions for forced response under periodic forcing is shown in Fig. \ref{SP_results}. By employing the oblique SSM reduction we can accurately track the backbone curve, and accurately predict the forced response curves of the system under harmonic forcing acting on both masses for several different forcing amplitudes. In contrast, the normal projection yields inaccurate results. 


\subsection{Shaw-Pierre system on an oscillating cart}
\begin{figure*}[b]
	\centering
	\includegraphics[scale=1.2]{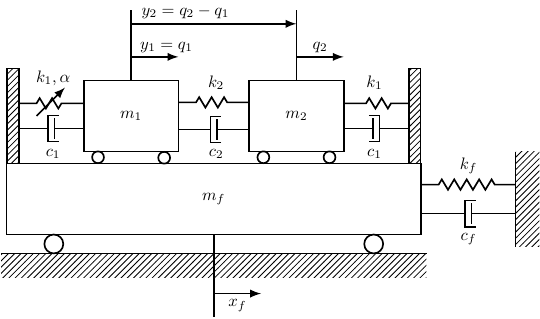}
	\caption{Geometry of the Shaw-Pierre mechanical system on a cart studied in \citet{haller_2024}. The coordinate $x_f$ refers to the absolute displacement of the cart, while $q_1$ and $q_2$ are the relative displacements of the two masses with respect to $x_f$.}
	\label{3dof_SP_system}
\end{figure*}

We adopt the three-degree-of-freedom model of \citet{haller_2024}, which consists of the Shaw-Pierre system we analyzed in \ref{SP_example}, but now mounted on an oscillating cart of mass $m_f$, as seen in Fig. \ref{3dof_SP_system}. 

The cart is connected to a wall via a linear spring and a linear damper. The equations of motion are still in the general form of \eqref{mechanical_system_equation}, but now with $\mathbf{y} = \left(y_1, y_2, x_f \right)^\mathrm{T}$ and
\begin{equation}
\begin{aligned}
&
\mathbf{M} = 
\begin{pmatrix}
m_1&0&m_1\\
m_2&m_2&m_2\\
0&0&m_f
\end{pmatrix},\quad
\mathbf{C} = 
\begin{pmatrix}
c_1 & -c_2 & 0 \\
c_1 & c_1 + c_2 & 0 \\
-2c_1 & -c_1 & c_f
\end{pmatrix}, \\ &
\mathbf{K} = 
\begin{pmatrix}
k_1 & -k_2 & 0 \\ 
k_1 & k_1 + k_2 & 0 \\ 
-2k_1 & -k_1 & k_f
\end{pmatrix},\quad
\mathbf{f}_\mathrm{nl} = 
\begin{pmatrix}
\alpha y_1^3 \\ 0 \\ -\alpha y_1^3
\end{pmatrix}.
\end{aligned}
\end{equation}
We fix the nondimensional parameter values 
\begin{equation*}
\begin{aligned}
    &
    m_1 = m_2 = m_f = 1, \quad c_1 = 0.05, \quad c_2 = c_f = 0.01, \\ 
    &k_1 = 1, \quad k_2 = 3.325, \quad k_3 = 33.25,\quad \alpha = 0.5,
\end{aligned}
\end{equation*}
for which the linearized system at the equilibrium has the eigenvalues $ \lambda_{1,2} = -0.022 \pm i\, 0.97, \quad \lambda_{3,4} = -0.035\pm i\,2.77, \quad \lambda_{5,6} = -0.059\pm i\,5.94$.
We want to construct a 2D reduced-order model on 2D primary SSM $\mathcal{W}\left(E\right)$ tangent to the slow spectral subspace $E$ corresponding to the eigenvalues $\lambda_{1,2}$.The real and the imaginary parts of the eigenvectors corresponding to $\lambda_{3,4}$ and  $\lambda_{5,6}$ span the 4D fast subspace $F$ in this example.

\begin{figure}[H]
    \centering
    \subfloat[\label{3dof_SP_norm_reconstruction}]{
        \includegraphics[scale=0.32]{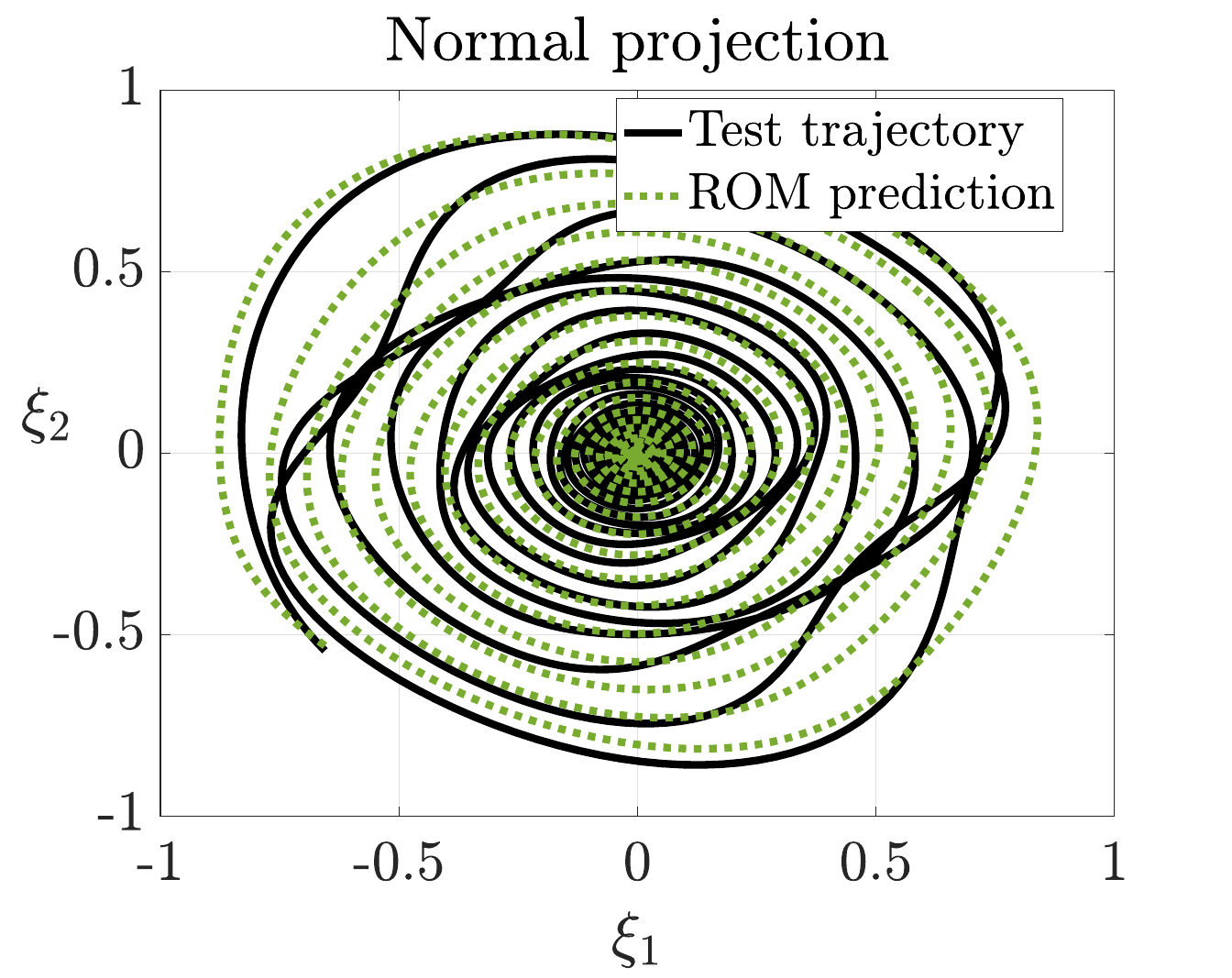}
    }
    \quad
    \subfloat[\label{3dof_SP_obl_reconstruction}]{
        \includegraphics[scale=0.32]{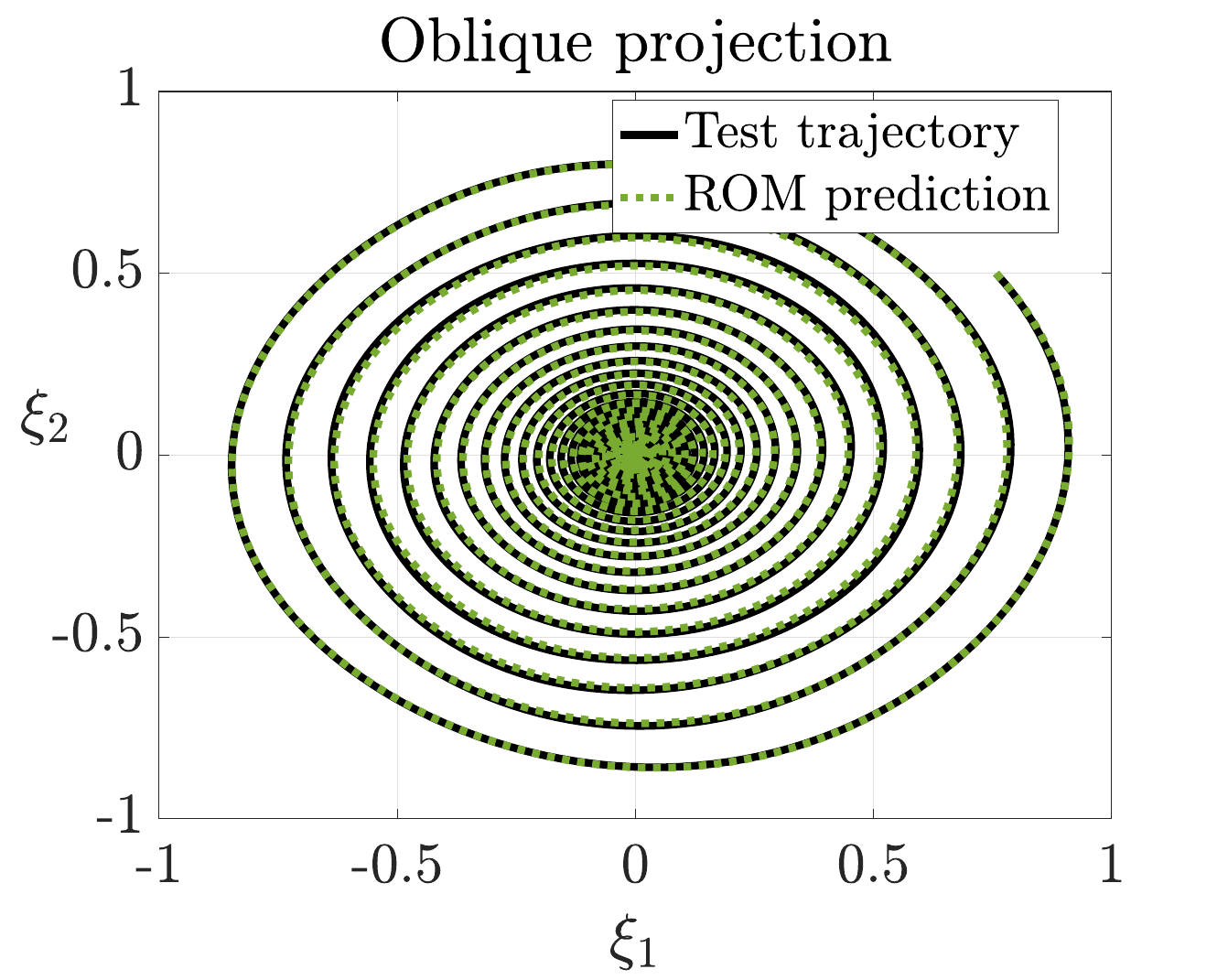}
    } \\
     \subfloat[\label{3dof_SP_norm_FRC}]{
        \includegraphics[scale=0.32]{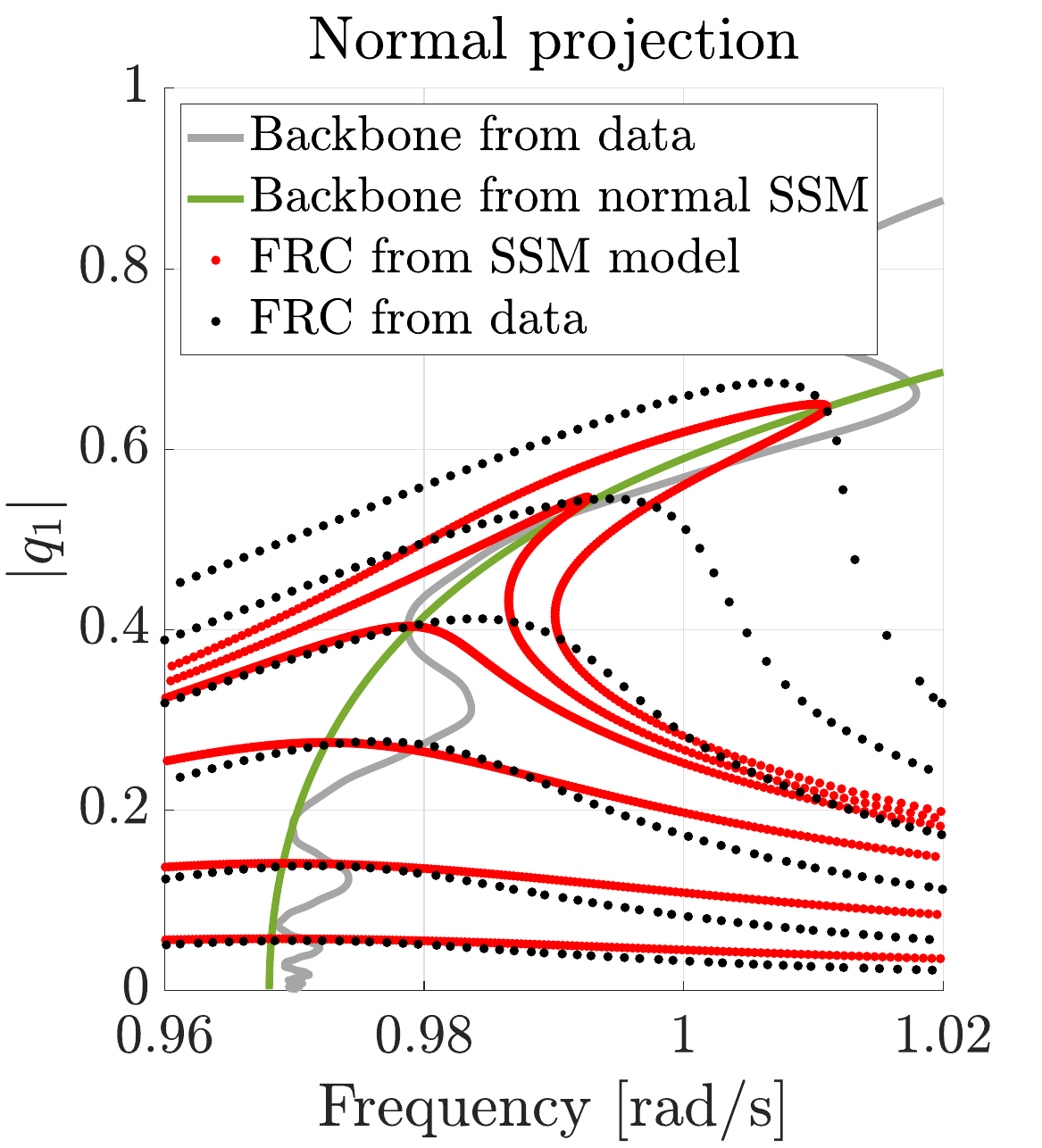}
    }
    \quad
    \subfloat[\label{3dof_SP_obl_FRC}]{
        \includegraphics[scale=0.32]{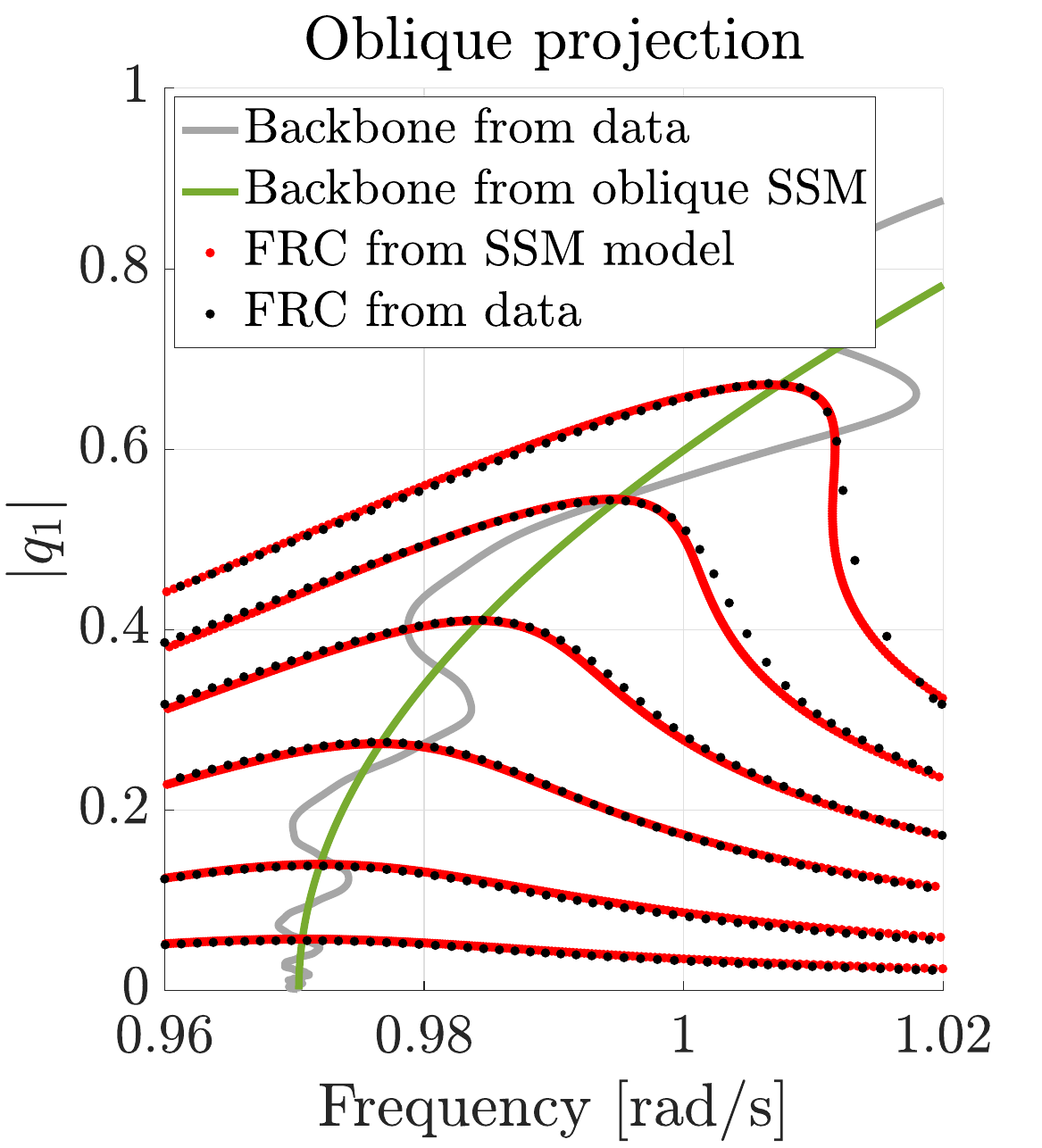}
    } 
    \caption{The reduced coordinates $\xi_1$ and $\xi_2$ are computed via normal and oblique projection onto the slow subspace of the system in Fig. \ref{3dof_SP_system}. The resulting SSM-reduced model (ROM) predictions are compared in terms of prediction of a test trajectory \subref{3dof_SP_norm_reconstruction} and \subref{3dof_SP_obl_reconstruction}, reconstruction of the backbone and prediction of the forced response curves \subref{3dof_SP_norm_FRC} and \subref{3dof_SP_obl_FRC}.} 
    \label{3dof_SP_results}
\end{figure}

We use SSMLearn to compute a $5^{\text{th}}$-order polynomial approximation of $\mathcal{W}\left(E\right)$ and its reduced dynamics from two trajectories observed in the full phase space (n = p = 6). The model's predictive ability is assessed on a test trajectory with an initial condition outside the training data, along with the reconstruction of the backbone curve and the prediction of the forced response curves under harmonic excitation of the three masses. The forcing mimics the shaking of the entire system, whereby all masses are subjected to the same forcing amplitude and phase, given $m_1 = m_2 = m_f$. As seen in Fig. \ref{3dof_SP_results}, the SSM constructed over obliquely projected reduced coordinates yields superior results compared to its normally projected counterpart. 

\subsection{Experimental data from nonlinear beam oscillations}
In our last example, we derive an SSM-reduced model from experimental data for the bending motion of a nonlinear beam. In this experiment, the response of the beam (transmitter) to external excitation induces gravitational forces on a bending beam resonator (receiver), located at a distance (\citet{brack2022}). The interaction between transmitter and receiver beams can then be studied to estimate the gravitational constant between the two objects. As such, accurate modeling of the beam becomes crucial. The forced response curves computed around the slowest bending modes of the transmitter beam reveal clear nonlinear behavior. Here we seek to construct a nonlinear 2D SSM-reduced model of the transmitter beam from decaying trajectory data.
\begin{figure*}[]
    \centering
    \subfloat[\label{gravitational_beam_setup}]{
        \includegraphics[scale=0.5]{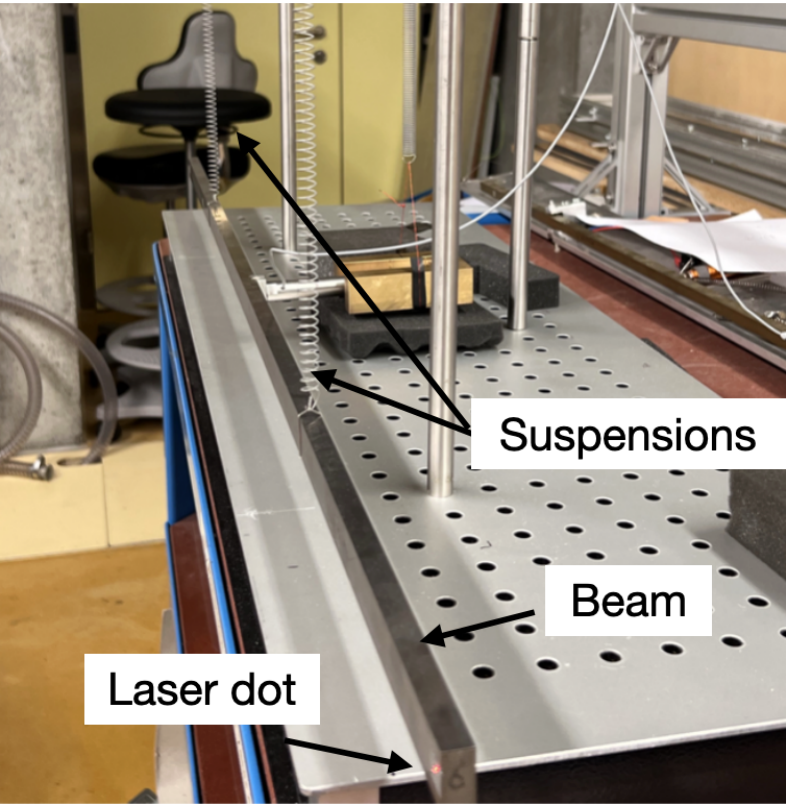}
    }
    \quad
    \subfloat[\label{grav_backbone}]{
        \includegraphics[scale=0.33]{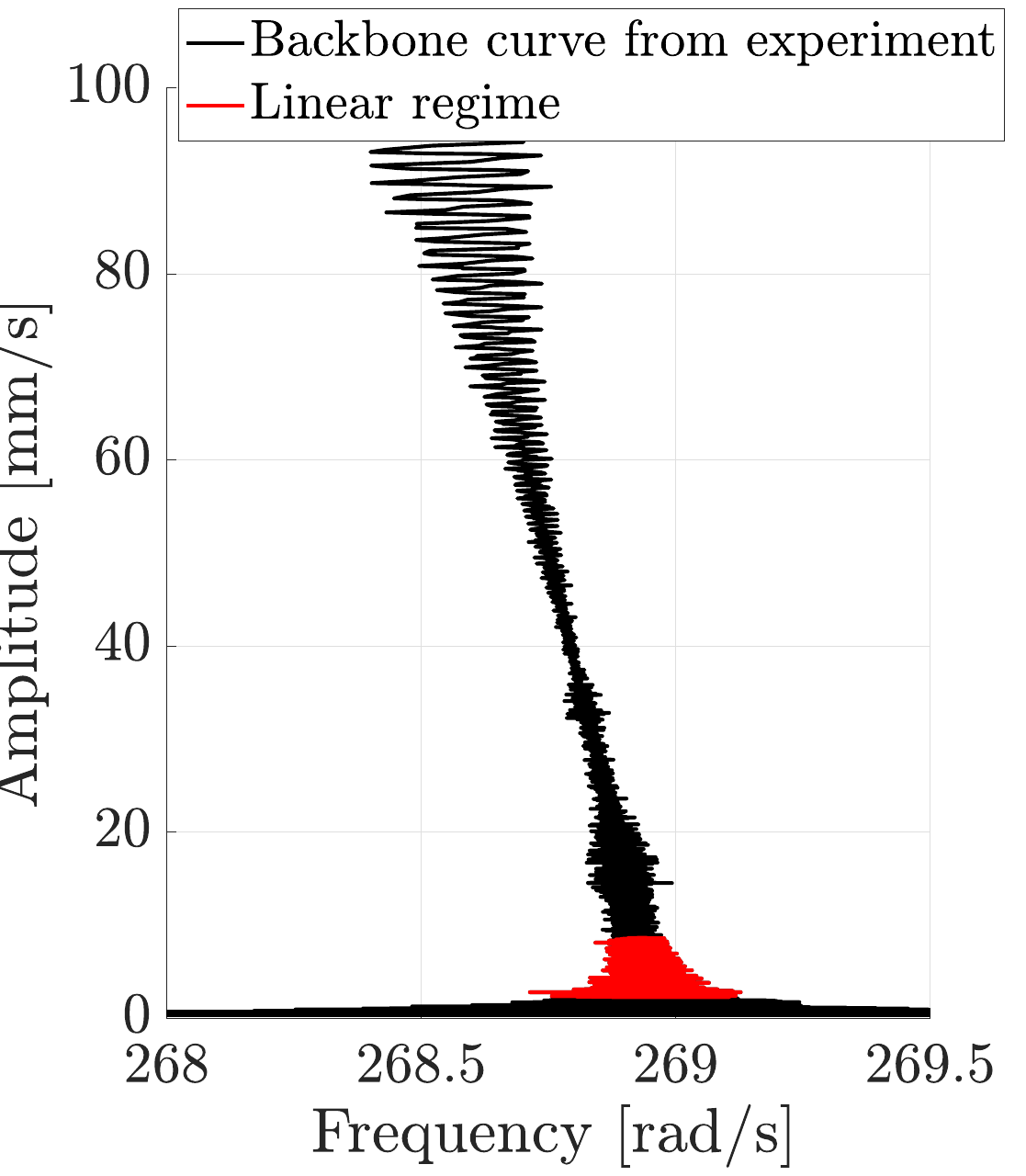}
    }
    \caption{Experimental setup and data of the vibrating nonlinear beam. \subref{gravitational_beam_setup} The experimental setup, with the transversal velocity measured at the tip of the beam via a laser. \subref{grav_backbone} The oscillatory backbone curve of a decaying trajectory, computed by means of the Peak Finding and Fitting (PFF) algorithm. The linear regime in which the oblique projection is computed is highlighted in red.} 
    \label{gravitational_experiment_data_results}
\end{figure*}

The tungsten beam has dimensions $1\,$m$\,\times\,20\,$mm$\,\times\,10\,$mm and a mass of $3,875.6\,$g. It is suspended by two strings attached at the nodal points of the first bending mode, as shown in Fig. \ref{gravitational_beam_setup}. A piezoelectric transducer is mounted at the middle point of the beam, where a counter mass amplifies the excitation. The velocity near one end of the beam is measured by a laser interferometer. Given a single scalar observable, we cannot satisfy Whitney's embedding theorem, for which we need $p > 2 d$ independent observables to embed a d-dimensional manifold. To address this, we invoke Takens' delay embedding theorem (\citet{takens_1981}), which states that delayed copies of a single scalar observables can serve as independent observables, thereby enabling the manifold embedding.

We construct the six-dimensional delay-embedded observable
\begin{equation}
\boldsymbol{\eta} = \left(v(t), v(t + \Delta t), \dots, v(t + (p-1)\Delta t) \right)^\mathrm{T}, \quad p = 6,
\end{equation}
where $v(t)$ denotes the single velocity measurement at the location of the laser dot and $\Delta t $ is the time delay. We took $\Delta t = 6 \cdot 10^{-4}$, which is the reciprocal of the acquisition frequency of the signal. The choice of $p$ satisfies the condition $p > 2d$ to embed the 2D slow SSM of dimension 2 via a generic observable (\citet{takens_1981}).

\begin{figure*}
    \centering
    \subfloat[\label{grav_FRC_norm}]{
        \includegraphics[scale=0.28]{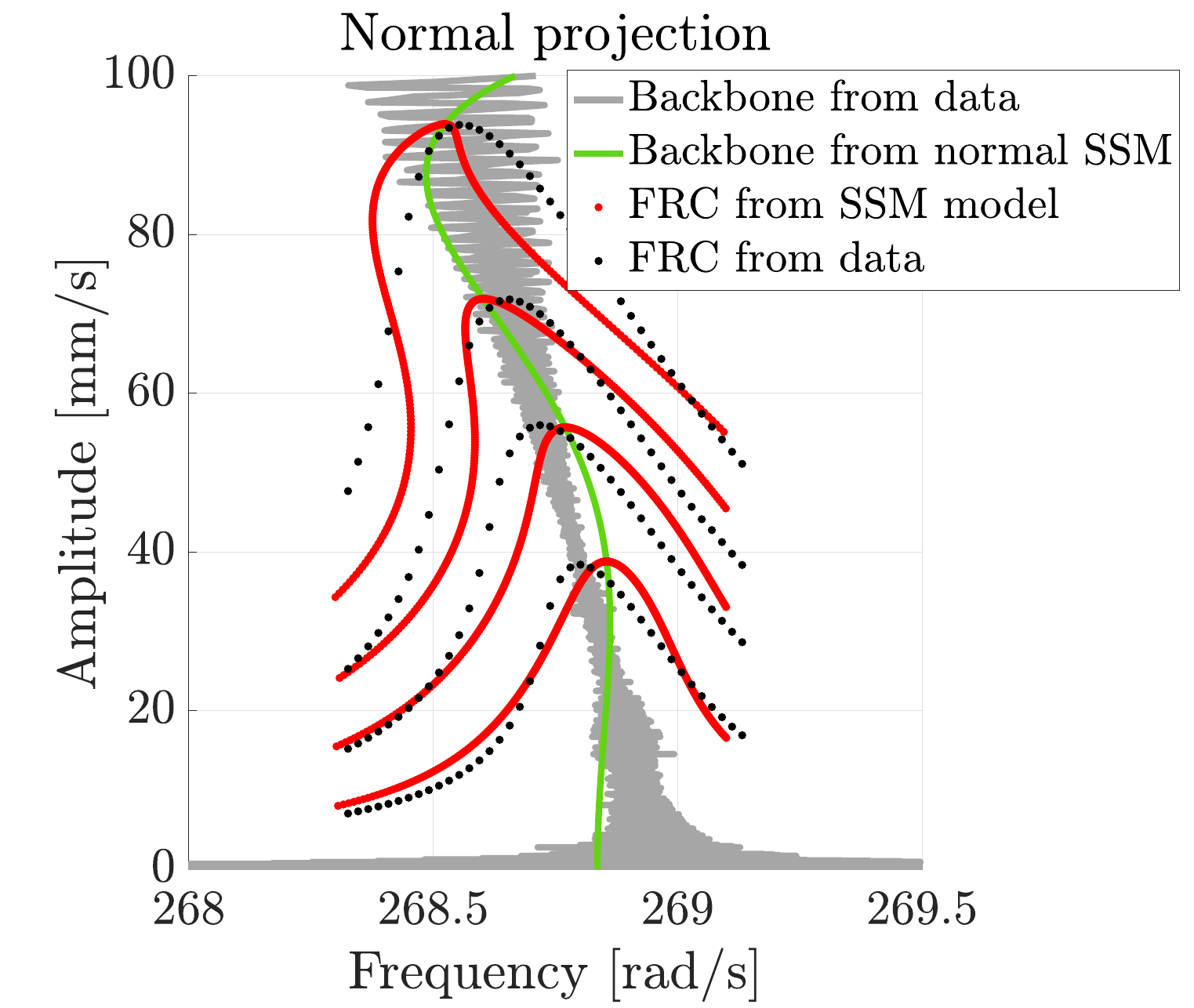}
    }
    \quad
    \subfloat[\label{grav_FRC_obl}]{
        \includegraphics[scale=0.28]{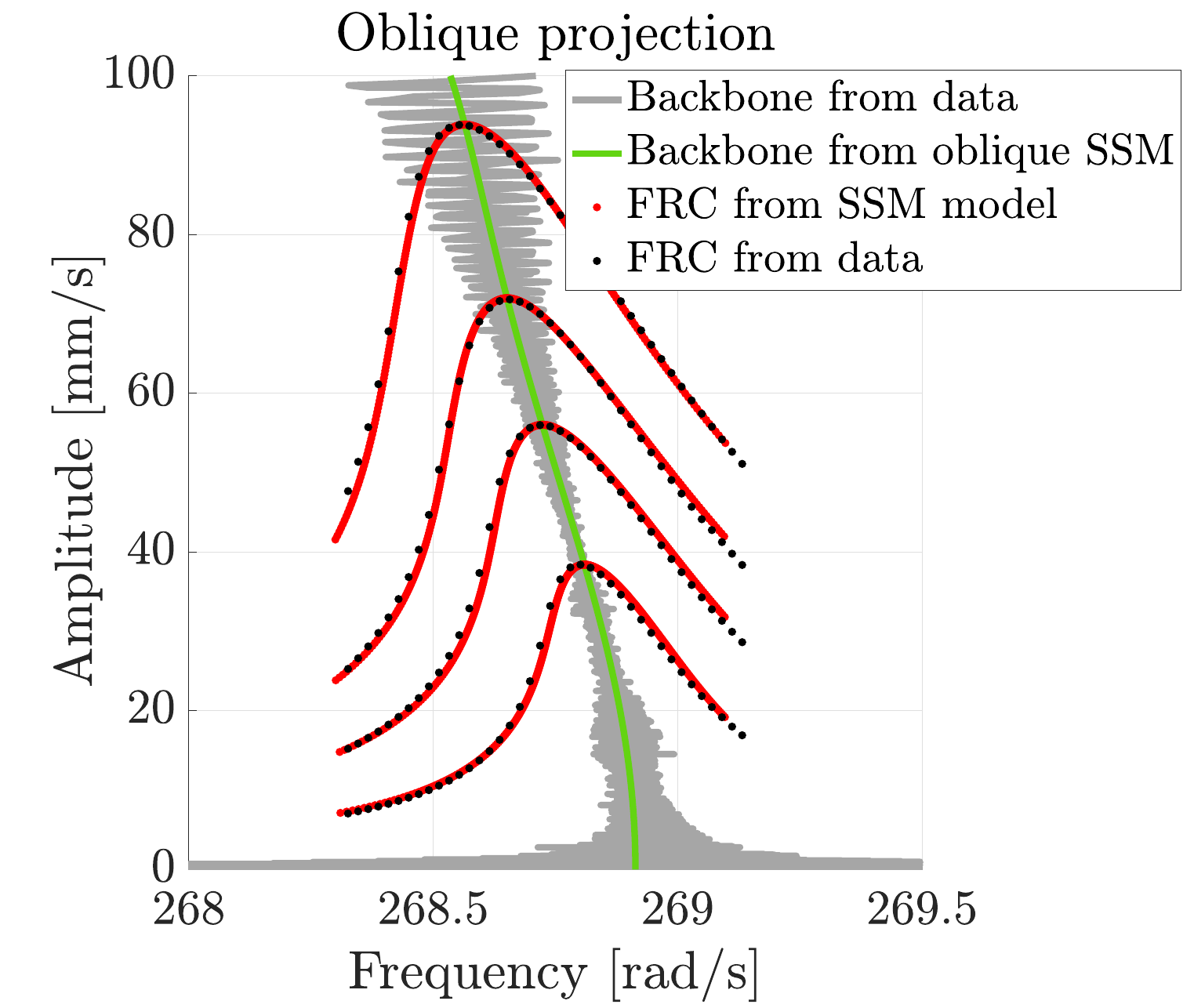}
    }\\
    \subfloat[\label{grav_phase_norm}]{
        \includegraphics[scale=0.28]{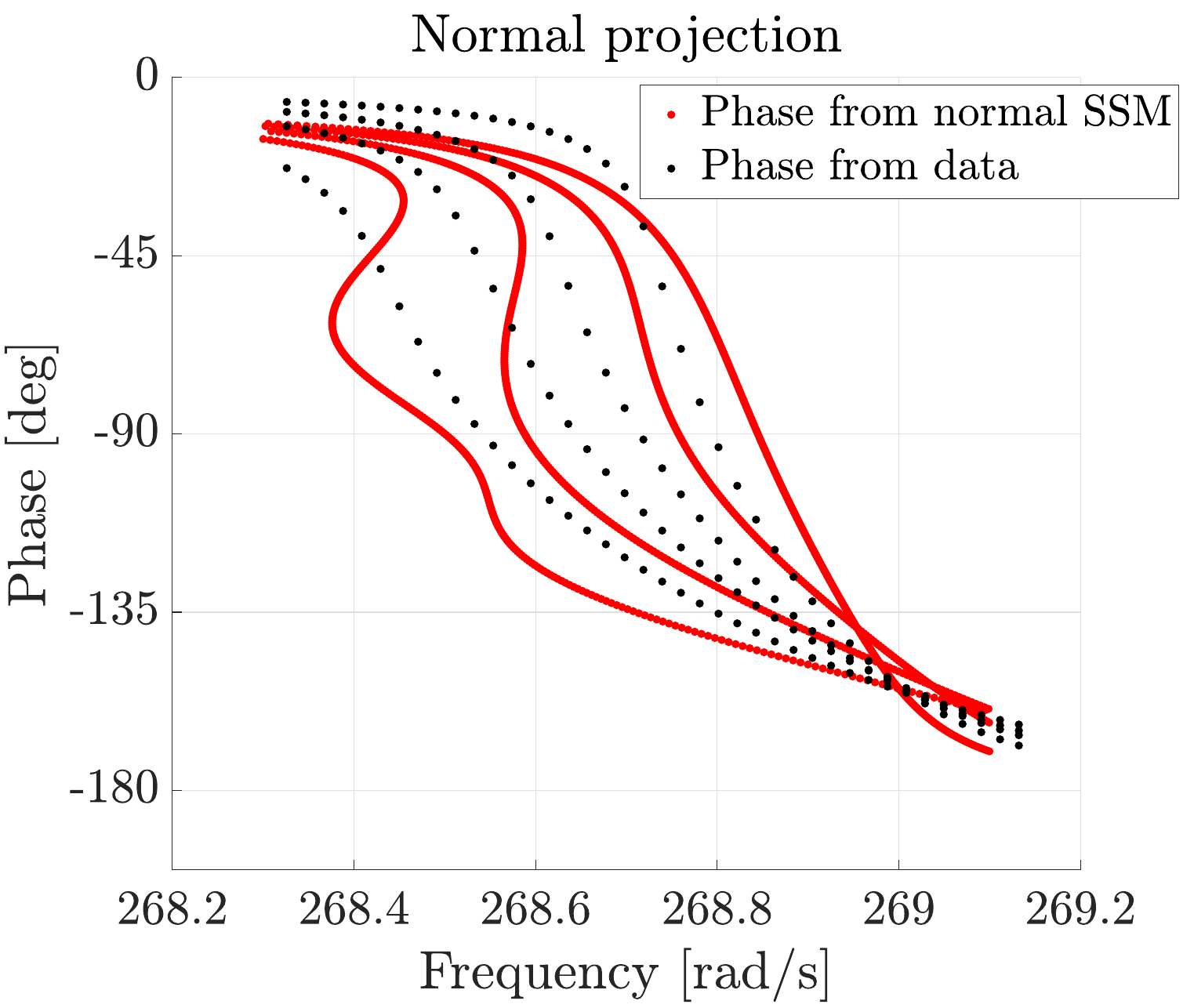}
    }
    \quad
    \subfloat[\label{grav_phase_obl}]{
        \includegraphics[scale=0.28]{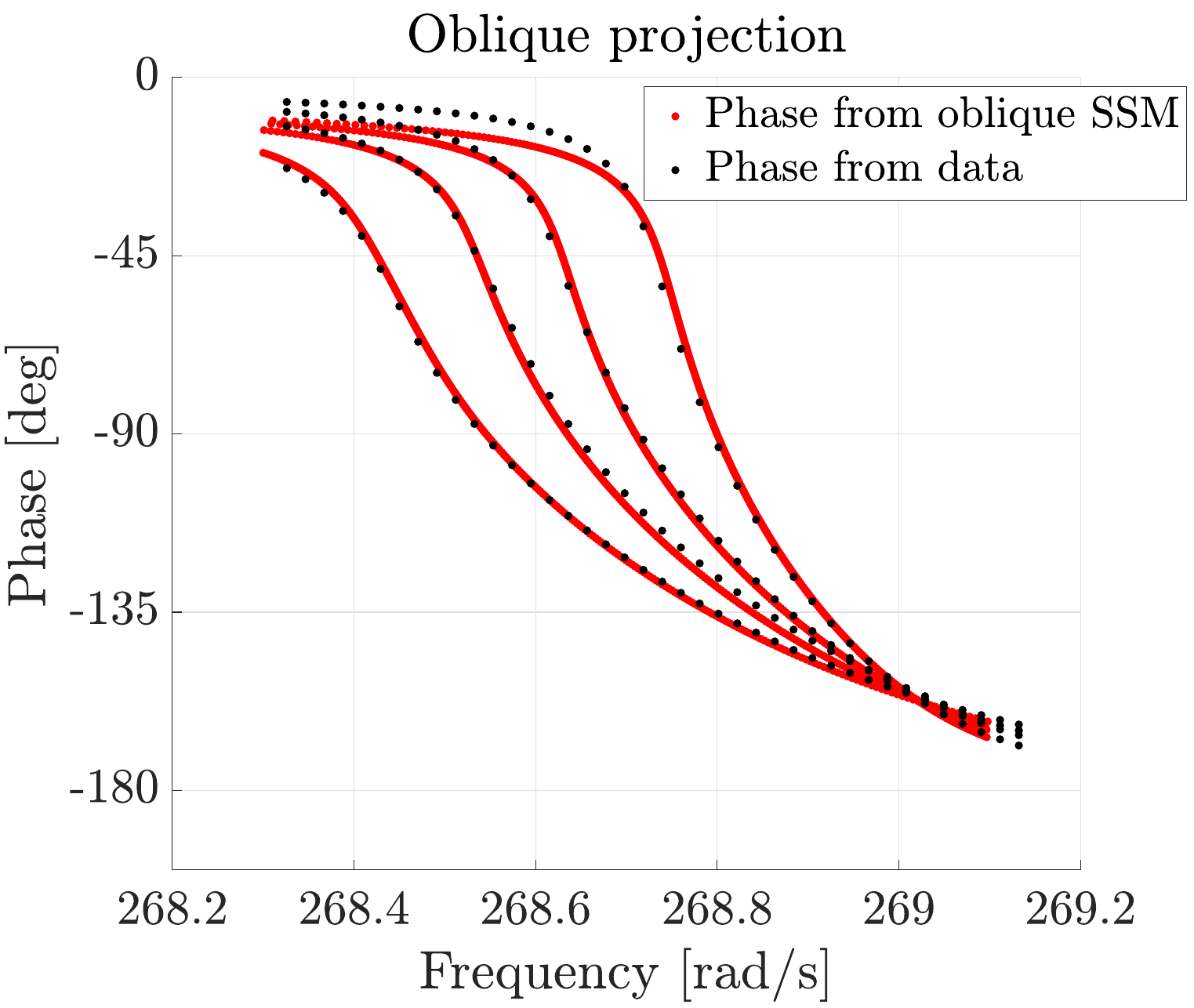}
    } 
    \caption{Comparison between the predictions of the primary SSM-based reduced-order model with normally and obliquely projected reduced coordinates, in terms of backbone curve and amplitude response (Subfigures \subref{grav_FRC_norm} and \subref{grav_FRC_obl}) and phase curves (Subfigures \subref{grav_phase_norm} and \subref{grav_phase_obl}).} 
    \label{gravitational_experiment_data_results}
\end{figure*}

We subject the beam to excitation at various forcing levels and in a range of frequencies around the first bending resonance. Recording the amplitude of the steady response as a function of the excitation frequency, we obtain the forced response curves shown in black dots in Figs. \ref{grav_FRC_norm} -\subref{grav_FRC_obl}. These curves represent our test set to be reconstructed by a 2D primary SSM-reduced model of the infinite-dimensional continuum beam. The training trajectory employed for the data-driven construction of this model consists of a single decaying trajectory obtained after a forcing at a specific amplitude and frequency is turned off. This trajectory will lie on a fractional SSM with probability one, as the full family of SSMs fills the entire domain of attraction of the unforced equilibrium. The primary SSM forms a single (measure zero) curve in this family.

We determine the backbone curve from the single decaying trajectory using the PFF algorithm. Shown in Fig. \ref{grav_backbone}, the backbone curve displays nonlinear softening behavior, and exhibits pronounced oscillations in frequency, even at low amplitudes.  This implies a significant non-normality for the linearized system based on our findings in Section \ref{motivation}.
As described in steps \ref{step_1}-\ref{step_4} of Section \ref{method}, we compute the oblique projection as the one that minimizes the oscillations of the backbone curve in the linear regime (shown in Fig. \ref{grav_backbone}). We show that this procedure only requires a single decaying trajectory. Indeed, given the oscillatory nature of the system, the trajectory explores the region of the phase space relevant for the SSM reconstruction and also provides information for the foliation approximation in terms of oscillatory backbone curve.

Both the parametrization of the primary SSM $\mathcal{W}\left(E\right)$ and its reduced dynamics are approximated by a $7^{\text{th}}$-order polynomial. As seen in Fig. \ref{gravitational_experiment_data_results}, the SSM-reduced model obtained via normal projection has difficulty in reconstructing the experimental backbone curve and predicting the forced response curves for amplitudes and phases. In contrast, the primary SSM-reduced model with oblique projection (trained on a single decaying trajectory) reproduces very accurately the experimental forced response and phase curves. The data have been measured via lock-in amplifier by introducing a fictitious $90^\circ$ phase shift with respect to a reference signal. The experimental phase curves in Fig. \ref{gravitational_experiment_data_results} have been plotted taking into account this shift.

\section{Conclusion}
We have enhanced the original SSMLearn algorithm for SSM-reduced modeling by employing an oblique projection onto slow spectral subspaces instead of the originally used normal projection. This results in a significant improvement in reconstructing backbone curves and predicting forced response curves for nonlinear oscillatory systems with substantial non-normality in their linear parts. Oblique SSM reduction is particularly relevant for experimental data, in which generic training trajectories lie on fractional SSMs.

We also give an explanation for the phenomenon of oscillating backbone curves, which turn out to arise even in linear systems. While the oscillations in the backbone curves are tempting to blame on noise or imperfections, the phenomenon in fact arises from an interplay between fractional SSMs and non-normality of slow and fast spectral subspaces.

Our methodology selects the reduced coordinates on the slow subspace via an oblique projection from the space of observables, along the span of the remaining fast spectral subspaces. Unlike prior approaches seeking to reconstruct the full fast subspace $F$, we only seek to identify the orthogonal complement $F^\perp$ of $F$ in the observable space. This enables a major reduction on the data needed in our approach, as $F^\perp$ has the same low dimension as the slow spectral subspace $E$. This allows us to use a single trajectory to construct the required oblique projection that minimizes the oscillations in the trajectory's backbone curve. We have applied this technique to various examples, both numerical and experimental, showing that the reduction on primary SSMs using obliquely projected reduced coordinates yields superior results compared to its normally projected counterpart.

The current method is constrained to linear oblique projection of the observable coordinates onto the slow subspace, over which the primary SSM is constructed. As such, the oblique projection is computed in the linearized regime, even for nonlinear dynamical systems, which results in a linear approximation of the stable foliation. While the leading factors influencing the difference between the dynamics on fractional and the primary SSMs appear to be related to the linear part, our future work will explore nonlinear oblique projections in order to enhance the accuracy of SSM-reduced modeling even further.

\section*{Conflict of Interest}
The authors have no conflicts to disclose.

\section*{Authors contributions}
\noindent GH designed the research. LB, BK and GH carried out the research. LB and BK developed the software and analyzed the examples. JD, BZ and LB designed the experiments. BZ and LB carried out the experiments. LB wrote the paper. GH, JD and BK reviewed the paper. GH led the research team.\\

\section*{Data availability}
The data discussed in this work are publicly available at \url{https://github.com/haller-group/SSMLearn}.

\section*{Code availability}
The code supporting the results of this work is publicly available at \url{https://github.com/haller-group/SSMLearn}.

\bibliography{bibliography}

\end{document}